%% file: Physics_Review.tex
\documentclass[12pt]{iopart}
\usepackage{iopams}
\usepackage{amssymb,amsfonts,graphicx}
\usepackage[nocompress]{cite}

\newcommand{\R}{\mathbb{R}}
\newcommand{\C}{\mathbb{C}}
\newcommand{\N}{\mathbb{N}}
\newcommand{\Z}{\mathbb{Z}}
\newcommand{\E}{\mathbb{E}}
\renewcommand{\H}{\mathcal{H}}
\newcommand{\F}{\mathcal{F}}
\newcommand{\K}{\mathcal{K}}
\newcommand{\J}{\mathcal{J}}
\newcommand{\ND}{N_{\mathrm{D}}}
\newcommand{\NN}{N_{\mathrm{N}}}

\newcommand{\eps}{\varepsilon}
\newcommand{\BigOh}{\mathcal{O}}
\renewcommand{\phi}{\varphi}
\renewcommand{\P}{\mathbb{P}}
\newcommand{\grad}{\mathop{\mathrm{grad}}}
\newcommand{\divg}{\mathop{\mathrm{div}}}
\newcommand{\Laplace}{\mathop{\triangle}}
\newcommand{\vol}{\mathop{\mathrm{vol}}}
\newcommand{\Res}{\mathop{\mathrm{Res}}\limits}
\newcommand{\dd}{\mathrm{d}}
\newcommand{\url}[1]{\texttt{#1}}
\newtheorem{lemma}{Lemma}[section]
\newtheorem{theorem}{Theorem}[section]
\newtheorem{corollary}{Corollary}[section]

\newtheorem{defn}{Definition}
\newtheorem{remark}{Remark}[section]

\renewcommand{\theenumi}{\roman{enumi}}

\begin{document}
\title[Laplace Operators on Fractals]{Laplace Operators on Fractals and Related
  Functional Equations}

\author{Gregory Derfel\textsuperscript{1}, Peter
  J. Grabner\textsuperscript{2}, Fritz Vogl\textsuperscript{3}}
\address{\textsuperscript{1}Department of Mathematics,
Ben Gurion University of the Negev,
Beer Sheva 84105,
Israel}

\address{\textsuperscript{2}
Institut f\"ur Analysis und Computational Number Theory,
Technische Universit\"at Graz,
Steyrergasse 30,
8010 Graz,
Austria}

\address{\textsuperscript{3}Institut f\"ur Analysis und Scientific Computing,
Technische Universit\"at Wien,
Wiedner Hauptstra\ss{}e 8--10,
1040 Wien,
Austria}

\eads{\mailto{derfel@math.bgu.ac.il}, \mailto{peter.grabner@tugraz.at},
\mailto{F.Vogl@gmx.at}}
\begin{abstract}
  We give an overview over the application of functional equations,
  namely the classical Poincar\'e and renewal equations, to the study
  of the spectrum of Laplace operators on self-similar fractals. We
  compare the techniques used to those used in the euclidean
  situation. Furthermore, we use the obtained information on the
  spectral zeta function to define the Casimir energy of fractals. We
  give numerical values for this energy for the Sierpi\'nski gasket.
\end{abstract}

\submitto{\JPA}

\maketitle

\tableofcontents
\title[Laplace Operators on Fractals]{}

\input{Introduction}
\input{Fractals}
\input{Laplacian}
\input{Spectral}
\input{Self-similarity}

\ack
The first author is partially supported by the Israel Science Foundation (ISF).
The second author is supported by the Austrian Science Fund project S9605.

The authors would like to thank D.~Berend, S.~Molchanov, and M.~Solomyak
for valuable remarks and helpful discussions.
\section*{References}
\addcontentsline{toc}{section}{References}
\bibliographystyle{unsrt} 
\bibliography{review}
\end{document}

%% file: Introduction.tex
\section{Introduction}\label{sec:intr}

\subsection{Historical perspective}\label{sec:historical-perspective}
Strange objects, which are poorly characterised by their topological
dimension, such as the Weierstrass continuous, but nowhere
differentiable function, the van Koch curve, the Sierpi\'nski gasket
(see Figure~\ref{fig:gasket}), the Sierpi\'nski carpet (see
Figure~\ref{fig:carpet}), etc., have been familiar in mathematics for
a long time.  Now such objects are called {\em fractals}.

The word ``fractal'' was coined by Mandelbrot in the 1970s.  His
foundational treatise \cite{Mandelbrot1982:fractal_geometry_nature}
contains a great diversity of examples from mathematics and natural
sciences. Coastlines, topographical surfaces, turbulence in fluid, and
so on are just few instances from the variety of natural objects that
may be described as fractals. There is no generally agreed exact
definition of the word ``fractal''.

The fractals studied in the context of analysis on fractals are all
\emph{self-similar}.  Moreover, we deal here mainly with
\emph{finitely ramified} fractals, i.e.,  fractals that can be
disconnected by removing a specific finite number of points.

Initial interest in processes (analysis) on fractals came from
physicists working in the theory of {\em disordered media}.  It turns
out that heat and wave transfer in disordered media (such as polymers,
fractured and porous rocks, amorphous semiconductors, etc.) can be
adequately modelled by means of fractals and random walks on them. (See
the initial papers by Alexander and Orbach
\cite{Alexander_Orbach1982:density_states_fractals}, Rammal and
Toulouse \cite{Rammal_Toulouse1983:random_walks_fractal}.  See also
the survey by Havlin and Ben-Avraham
\cite{Havlin_Ben-Avraham1987:diffusion_disordered_media} and the book
by the same authors
\cite{Havlin_Ben-Avraham2000:diffusion_reactions_fractals} -- for an
overview of the now very substantial physics literature and
bibliography.)

Motivated by these works,  mathematicians got interested in developing
the of ``analysis on fractals''. For instance, in order to analyse how
heat diffuses in a material with fractal structure, one needs to
define a ``heat equation'' and a ``Laplacian'' on a fractal. The
problem contains somewhat contradictory factors. Indeed, fractals like
the Sierpi\'nski gasket, or the van Koch curve do not have any smooth
structures and one cannot define differential operators on them
directly. We will comment on the difference between the euclidean and
the fractal situation throughout this paper.

\subsection{Probabilistic approach}\label{sec:probabilistic-approach}
In the mid 1980s probabilists solved the problem by constructing
``Brownian motion'' on the Sierpi\'nski gasket.  Goldstein
\cite{Goldstein1987:random_walks_diffusion}, Kusuoka \cite
{Kusuoka1987:diffusion_process_fractal}, and a bit later Barlow and
Perkins \cite{Barlow_Perkins1988:brownian_motion_sierpinski}, 
independently took the first step in the mathematical development of
the problem. Their method of construction is now called the {\em
  probabilistic approach}. Namely, they considered a sequence of
random walks $X^{(n)}$ on graphs $G_n$, which approximate the
Sierpi\'nski gasket $\Gamma$, and showed  that by taking a certain scaling
factor,  those random walks converge to a diffusion process on the
Sierpi\'nski gasket.

In order to obtain a nontrivial limit, the time scale for each step
should be the average time for the random walk on $G_{n+1}$, starting
from a point in $G_n$, to arrive  at a point in $G_n$,
except for the starting point. By the self-similarity and symmetry of
$\Gamma$ this average time $\tau$ is independent of $n$ and it is
equal to the average time for the random walk on $G_1$, starting from
$a$ to arrive at either $b$ or $c$. (Here $a$, $b$, $c$, are the
vertices of the initial equilateral triangle $G_0$.)  Elementary
calculations show that $\tau=5$.  According to
\cite{Goldstein1987:random_walks_diffusion,
  Kusuoka1987:diffusion_process_fractal,
  Barlow_Perkins1988:brownian_motion_sierpinski}  the processes
$2^{-n}X^{(n)}([5^nt])$ weakly converge (as $n \rightarrow\infty$) to a
non-trivial limit $X_t$ on $\Gamma$, which is called Brownian motion
on $\Gamma$. In this approach the Laplacian is defined as
an  \emph{infinitesimal generator} of $X_t$.


An important observation was made by Barlow and Perkins in \cite
{Barlow_Perkins1988:brownian_motion_sierpinski}.  Let $Z_n=T^{n.0}_1$
be the \emph{first hit} time by $X^{(n)}$ on $G_0$. Then $Z_n$ is a
\emph{simple branching process}.  Its off-spring distribution $\eta$
has the generating function $q(z)=\E(z^{\eta})=z^2/(4-3z)$ and, in
particular, $\E(\eta) = q'(1) = 5$.  Thus, $Z_n$ is a
\emph{super-critical} branching process.  It is known (see
\cite{Harris1963:theory_branching_processes}) that in this case
$5^{-n}Z_n$ tends to a limiting random variable $Z_\infty$.

The moment generating function of this random variable
\begin{equation*}
f(z)=\mathbb{E}e^{-zZ_\infty}
\end{equation*}
satisfies the functional equation
\begin{equation} \label{eq:branching} f(\lambda z)=q(f(z)),
\end{equation}
which is the \emph{Poincar\'e equation} (see also
Section~\ref{sec:expl-self-simil}  below).

The Poincar\'e equation associated with the Brownian motion turns out
to be a very useful tool in the study of the  detailed properties (for example,
heat kernel) of the Brownian motion.

Lindstr\o m \cite{Lindstroem1990:brownian_motion_nested} extended the
construction of the Brownian motion from the Sierpi\'nski gasket to more
general \emph{nested fractals}, (which are finitely ramified
self-similar fractals with strong symmetry).  The Lindstr\o m
\emph{snowflake} is a typical example of  a nested fractal (see
Figure~\ref{fig:snowflake}).

Readers may refer to Barlow's lecture notes
\cite{Barlow1998:diffusion_fractals } for a self-contained survey of
the probabilistic approach.
\subsection{Anomalous diffusion}\label{sec: anomalous-diffusion}

It has been discovered in an early stage already (see
\cite{Alexander_Orbach1982:density_states_fractals,
  Rammal_Toulouse1983:random_walks_fractal,
  Goldstein1987:random_walks_diffusion,
  Kusuoka1987:diffusion_process_fractal,
  Barlow_Perkins1988:brownian_motion_sierpinski}) that diffusion on
fractals is \emph{anomalous}, different than that in a regular
space. For a regular diffusion, or (equivalently) a simple random walk in
all integer dimensions $d$, mean-square displacement is proportional to
the number of steps $n$: $\E(X_n)^2 =c n$ (Fick's law, 1855).  On the
other hand, in the case of the Sierpi\'nski gasket $\E^{x}(X_n -x)^2
\asymp n^{2/\beta} $, where $\asymp$ means that the ratio between the two sides 
is `bounded above and below by
positive constants' and $\beta =\lg5/\lg2$ is called the \emph{walk
  dimension}.

This slowing down of the diffusion is caused, roughly speaking, by the
removal of large parts of the space.

De Gennes \cite{deGennes1976:percolation_concept} was amongst the
first, who realised the broad importance of anomalous diffusion and
coined the suggestive term ``the ant in the labyrinth'', describing
the meandering of a random walker in percolation clusters.
\subsection{Analytic approach}\label{sec:analytic-approach}

The second approach, based on difference operators, is due to Kigami
\cite{Kigami1989:harmonic_calculus_sierpinski}.  Instead of the
sequence of random walks, one can consider a sequence of discrete
Laplacians on a sequence of graphs, approximating the fractal.  It is
possible to prove that under a proper scaling these discrete
Laplacians would converge to an  ``well-behaved'' operator with dense domain,
called the
Laplacian on the Sierpi\'nski gasket. This alternative approach is
usually called the \emph{analytic approach}.

Later it was extended by Kigami
\cite{Kigami1993:harmonic_calculus_pcf,
  Kigami1998:distributions_localized_eigenvalues,
  Kigami2001:analysis_fractals} to more general class of fractals --
\emph{post critically finite self-similar sets} (p.c.f), which
roughly correspond to finitely ramified self-similar fractals.

The two approaches described above are complementary to each other.

The advantage of the analytic approach is that one gets concrete and
direct description of harmonic functions, Laplacians, Dirichlet forms,
etc.  (See also
\cite{Fukushima_Shima1992:spectral_analysis_sierpinski,
  Kusuoka1989:dirichlet_forms_random_matrices}.)

On the other hand, however, the probabilistic approach is better
suited for the  study of  heat kernels.  Moreover, this approach can be
applied to \emph{infinitely ramified} self-similar fractals, which
include the Sierpi\'nski carpet, as a typical example (cf., 
\cite {Barlow_Bass1989:construction_brownian_motion}).
  
\subsection{Functional equations in the analysis on fractals}
\label{sec:functional-equations}
One important example, where the Poincar\'e equation arises in
connection with analysis on fractals, has been mentioned  already at
the end of Section~\ref{sec:probabilistic-approach}.

Further applications of functional equations in this field are related
to spectral zeta function $\zeta_\Delta$ of  the Sierpi\'nski gasket
and other more general fractals having \emph{spectral decimation}.
The phenomenon of spectral decimation was first observed and studied
by Fukushima and Shima
(\cite{Fukushima_Shima1992:spectral_analysis_sierpinski,
  Shima1993:eigenvalue_problem_laplacian,
  Shima1996:eigenvalue_problems_laplacians}) and further progress has
been made by Malozemov and
Teplyaev~\cite{Malozemov_Teplyaev2003:self_similarity_operators} and
Strichartz~\cite{Strichartz2003:fractafolds_spectra}.

The definition of spectral decimation is given in
Section~\ref{sec:spectral-decimation} (see Definition ~\ref{def1} below).
It implies, in particular, that eigenvalues of the Laplacian $\Delta$
on the fractal, which admits spectral decimation, can be calculated by
means of a certain polynomial $p(z)$, or rational function $R(z)$.
Hence, the spectral zeta function $\zeta_\Delta$ may be defined by means
of iterations $p^{(n)}$ of $p$, or $R^{(n)}$ of $R$ (see
~\cite{Teplyaev2004:spectral_zeta_function}).

The above mentioned iteration process, as is well known in 
iteration theory, may be conveniently described by the corresponding
Poincar\'e equation:
\begin{equation}  \label{eq:decimation_iteration}
\Phi(\lambda z)=p(\Phi(z)),
\end{equation}
where $\lambda=p'(0)>1$ (see, for example,
Beardon~\cite{Beardon1991:iteration_rational_functions} or
Milnor~\cite{Milnor2006:dynamics_complex}).

Using that, in Section~\ref{sec:spectr-zeta-funct}, we obtain the
meromorphic continuation of the zeta function $\zeta_\Delta$ into the
whole complex plane on the basis of the knowledge of the asymptotic behaviour
of the Poincar\'e function $\Phi$ in certain angular regions.

The poles of the spectral zeta function are called the \emph{spectral
  dimensions} (see ~\cite{Lapidus1992:spectral_fractal_geometry,
  Lapidus_Frankenhuysen2000:fractal_geometry_number,
  Lapidus_Frankenhuijsen2006:fractal_geometry_complex}).  For the physical
consequences of complex dimensions of fractals -- see (\cite{
  Akkermans_Dunne_Teplyaev2009:complex_dimensions,
  Akkermans_Dunne_Teplyaev2010:thermodynamics_fractals}).
In Section~\ref{sec:casimir}, we use the Poincar\'e function $\Phi$
for  the calculation of \emph{Casimir energy} on a fractal.

Finally, one can expect that functional equations with rescaling
naturally come about from problems, where renormalisation type
arguments are used to study self-similarity.  Furthermore, functional
equations, in contrast with differential ones, do not require any
\emph{smoothness} of solutions; they may possess nowhere differential
solutions, for example.

\subsection{Notes and remarks}\label{sec:notes}
There is now a number of excellent books, lecture notes and surveys on
different aspects of analysis and probability on fractals
\cite{Barlow1998:diffusion_fractals, Kigami2001:analysis_fractals,
  Strichartz2006:differential_equations_fractals,
  Lapidus_Frankenhuijsen2006:fractal_geometry_complex,
  Triebel1997:fractals_and_spectra, Kirillov2009:fractals,
  Kumagai2008:recent_developments_fractals,
  Teplyaev2010:diffusions_and_spectral,
  Havlin_Ben-Avraham2000:diffusion_reactions_fractals}.

It is next to impossible to describe all activities in this area.  The
objective of the present paper is different.  We restrict ourselves to
a brief overview of various approaches in the study of the Laplacian and its
spectral properties on certain self-similar fractals, and we put much
emphasis on the deep connection between the latter problem and functional
equations with rescaling, and the classical Poincar\'e equation, in
particular.

The Poincar\'e equation plays a very important role in the
mathematical theory of dynamical systems and in  iteration theory,
in particular, but is still much less known in the physics
literature. Hopefully, the current presentation of this realm of
problems intended for a general audience may fill this gap.

We begin our presentation of the Poincar\'e equation in
Section~\ref{sec:expl-self-simil} with a brief description of the
general case, when the coefficients of the polynomial $P$ and  the scaling factor
$\lambda$ are \emph{complex} numbers, and only afterwards turn to a
detailed discussion of the \emph{real} case. So far, the real case
only has arisen in analysis on fractals.  However, we think that the
general case is interesting on its own, and, probably, will also find
applications in the future.

In this overview we do not or only cursorily touch the following important topics:
\begin{itemize}
\item analysis on infinitely ramified fractals such as the Sierpi\'ski
  carpet \cite{Barlow_Bass1989:construction_brownian_motion,
    Barlow_Bass1999:brownian_sierpinski_carpet}. The very recent
  progress that has been made in proving the uniqueness of the
  diffusion on the Sierpi\'nski carpets
  \cite{Barlow_Bass_Kumagai+2010:uniqueness_of_brownian} provided a
  unification of the different approaches to diffusion on this class
  of fractals.
\item heat kernel long time behaviour and Harnack inequalities on
  general underlying spaces, including Riemannian manifolds, graphs
  and fractals, as special cases
  \cite{Barlow_Coulhon_Grigoryan2001:manifolds_graphs_with,
    Grigoryan2010:heat_kernels_metric,
    Grigoryan_Hu_Lau2010:comparison_inequalities,
    Barlow_Grigoryan2009:heat_kernel_upper}
\item analysis by means of potential theory and functional spaces (Sobolev
  or Besov spaces) techniques
  \cite{Naimark_Solomyak1994:eigenvalue_behaviour,
    Naimark_Solomyak1995:eigenvalue_behaviour,
    Naimark_Solomyak2001:eigenvalue_distribution_fractal,
    Solomyak_Verbitsky1995:spectral_problem_related,
    Triebel1997:fractals_and_spectra,
    Triebel2008:fractal_analysis_approach,
    Hu_Lau_Ngai2006:laplace_operators_related}

 We refer the reader to the original papers.
\end{itemize}

%% file: Fractals.tex
\section{Fractals and iterated function systems}
\label{sec:fractals}
Let us first recall that for any finite set of linear contractions on $\R^d$
\begin{equation*}
F_i(\mathbf{x})=\mathbf{b}_i+\mathbf{A}_i(\mathbf{x}-\mathbf{b}_i),\quad 
i=1,\ldots,m
\end{equation*}
with fixed points $\mathbf{b}_i$ and contraction matrices $\mathbf{A}_i$
($\|\mathbf{A}_i\|<1$, $i=1,\ldots,m$) there exists a unique compact set
$K\subset\R^d$ satisfying
\begin{equation}\label{eq:fixed-K}
K=\bigcup_{i=1}^m F_i(K).
\end{equation}
In general, the set $K$ obtained as fixed point in \eref{eq:fixed-K} is a
``fractal'' set
(cf.~\cite{Falconer1986:geometry_fractal_sets,
Hutchinson1981:fractals_self_similarity}). If the matrices $\mathbf{A}_i$ are
only assumed to be contracting, finding the Hausdorff-dimension
of $K$ is a rather delicate
(cf.~\cite{Falconer1988:hausdorff_dimension_self,
Falconer1992:dimension_self_affine_ii}) problem. 

For this paper we will make the following additional assumptions on the family
of contractions $F=\{F_i\mid i=1,\ldots,m\}$:
\begin{enumerate}
\item\label{simil} the matrices $\mathbf{A}_i$ are similitudes with factors
  $\alpha_i<1$: 
  \begin{equation*}
    \forall\mathbf{x}\in\R^d:\|\mathbf{A}_i\mathbf{x}\|=\alpha_i\|\mathbf{x}\|
  \end{equation*}
\item\label{open-set} $F$ satisfies the open-set-condition
  (cf.~\cite{Falconer1986:geometry_fractal_sets,
    Falconer1997:techniques_in_fractal}), namely, there exists a bounded open
  set $\mathcal{O}$ such that
  \begin{equation*}
    \bigcup_{i=1}^m F_i(\mathcal{O})\subset\mathcal{O}
  \end{equation*}
  with the union disjoint.
\end{enumerate}

Assuming that the maps in $F$ are similitudes and satisfy the
open-set-condition, the Hausdorff-dimension of the compact set $K$ given by
\eref{eq:fixed-K} equals the unique positive solution $s=\rho$ of the equation
\begin{equation}\label{eq:dimension}
\sum_{i=1}^m\alpha_i^s=1.
\end{equation}
The projection map
\begin{equation}\label{eq:pi}
  \begin{array}{lll}
    \pi:&\{1,\ldots,m\}^{\N}&\to K\\
    &(\eps_1,\eps_2,\ldots)&\mapsto \lim_{n\to\infty}F_{\eps_1}\circ
    F_{\eps_2}\circ\cdots\circ F_{\eps_n}(\mathbf{x})
  \end{array}
\end{equation}
(the limit is independent of $\mathbf{x}\in\R^d$) defines a ``parametrisation''
of $K$. If $\{1,\ldots,m\}^{\N}$ is endowed with the infinite product measure
$\mu$ given by
\begin{equation}\label{eq:mu}
\mu\left(\left\{(\eps_1,\eps_2,\ldots)\mid \eps_1=i_1,\eps_2=i_2,\ldots
\eps_k=i_k\right\}\right)=
\left(\alpha_{i_1}\alpha_{i_2}\cdots\alpha_{i_k}\right)^\rho,
\end{equation}
then $\pi$ is a $\mu$-almost sure bijection. The normalised $\rho$-dimensional
Hausdorff-measure on $K$ is then given by
\begin{equation}\label{eq:pi*mu}
\H_K^\rho(E):=
\frac{\H^\rho(E\cap K)}{\H^\rho(K)}=\pi_*(\mu)(E)=\mu\left(\pi^{(-1)}(E)\right).
\end{equation}
This is the unique normalised measure satisfying
(cf.~\cite[Theorem~28]{Falconer1997:techniques_in_fractal})
\begin{equation}\label{eq:similar-measure}
\H_K^\rho(E)=\sum_{i=1}^m\alpha_i^\rho\H_K^\rho(F_i^{(-1)}(E)).
\end{equation}

By the above discussion the functions
\begin{equation*}
\eps_n(\mathbf{x})=\left(\pi^{(-1)}(\mathbf{x})\right)_n
\end{equation*}
are defined $\H_K^\rho$-almost everywhere on $K$. We set
\begin{equation}\label{eq:eps_n}
\eps_n(\mathbf{x})=\min\left(\pi^{(-1)}(\{\mathbf{x}\})\right)_n
\end{equation}
to define them everywhere. The set $K$ together with the address space
$\Sigma=\{1,\ldots,m\}^{\N}$ and the maps $F_i$ ($i=1,\ldots,m$) defines a
\emph{self-similar structure} $(K,\Sigma,(F_i)_{i=1}^m)$.

In order to allow sensible analysis on the fractal set $K$, we need some
further properties. Especially, since we will later study diffusion processes,
we need $K$ to be connected. On the other hand, the techniques introduced later
require a finite ramification property usually called post-critical finiteness
(p.~c.~f.).

Connectivity of $K$ is characterised by the fact that for any pair $(i,j)$,
there exist $i_1,i_2,\ldots,i_n\in\{1,\ldots,m\}$ with $i=i_1$ and $j=i_n$ such
that $F_{i_\ell}(K)\cap F_{i_{\ell+1}}(K)\neq\emptyset$ for $\ell=1,\ldots,n-1$
(cf.~\cite{Hata1985:structure_self-similar}).

\begin{defn}\label{def:pcf}
Let $(K,\Sigma,(F_i)_{i=1}^m)$ be a self-similar structure. Then the set
\begin{equation*}
C=\pi^{-1}\left(\bigcup_{i\neq j}F_i(K)\cap F_j(K)\right)
\end{equation*}
is called the \emph{critical set} of $K$. The \emph{post-critical set} of $K$
is defined by
\begin{equation*}
P=\bigcup_{n=1}^\infty \sigma^n(C),
\end{equation*}
where $\sigma:\Sigma\to\Sigma$ denotes the shift map on the address space
$\Sigma$. If $P$ is a finite set, then $(K,\Sigma,(F_i)_{i=1}^m)$ is called
\emph{post-critically finite} (p.~c.~f.). This is equivalent to the finiteness
of $C$ together with the fact that all points of $C$ are ultimately periodic.
\end{defn}

The following sequence $V_m$ of finite sets will be used in
Section~\ref{sec:lapl-oper-dirichlet} to define a sequence of electrical
networks giving a harmonic structure on $K$. For more details we refer to 
\cite[Chapter~1]{Kigami2001:analysis_fractals}
\begin{defn}\label{def:vm}
Let $(K,\Sigma,(F_i)_{i=1}^m)$ be a post-critically finite self-similar
structure and $P$ its post-critical set. Let $V_0=\pi(P)$ and define $V_m$
iteratively by
\begin{equation*}
V_{n+1}=\bigcup_{i=1}^mF_i(V_n).
\end{equation*}
The sets $V_n$ are then finite, increasing ($V_n\subset V_{n+1}$) and
\begin{equation*}
K=\overline{\bigcup_{n\geq0}V_n}.
\end{equation*}

\begin{figure}[h]
  \centering
  \includegraphics[width=0.5\textwidth]{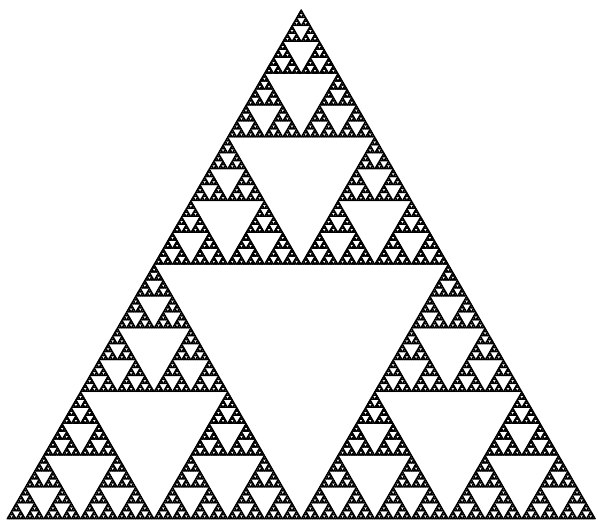}
  \caption{The Sierpi\'nski gasket; the points of $V_0$ are the vertices
    of the triangle}
  \label{fig:gasket}
\end{figure}

\begin{figure}[h]
  \centering
  \includegraphics[width=0.8\textwidth]{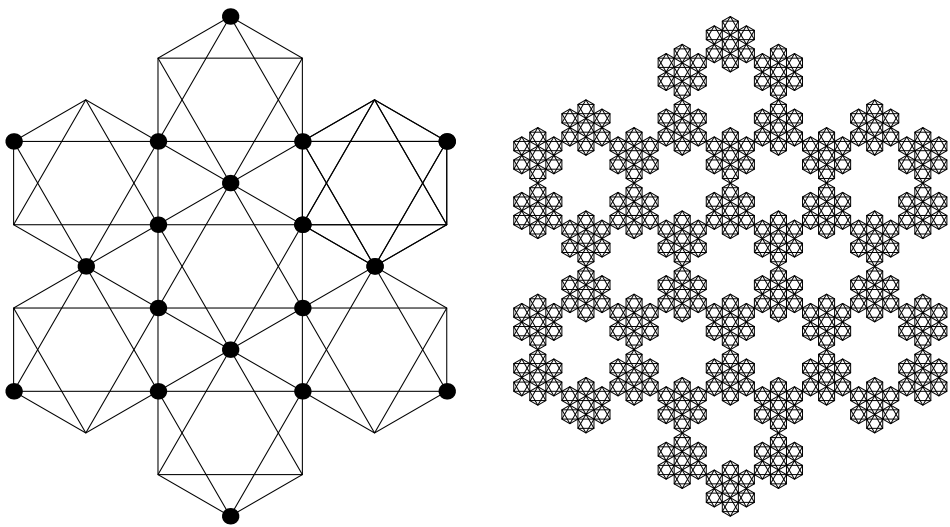}
  \caption{The Lindstr\o m snowflake with the corresponding set $V_0$}
  \label{fig:snowflake}
\end{figure}

\end{defn}

%% file: Laplacian.tex
\section{Laplace operators on fractals}\label{sec:laplace}
\subsection{Laplace operators on compact manifolds}
\label{sec:lapl-oper-comp}
Before we introduce the Laplace operator on certain classes of
self-similar fractals, let us shortly discuss the situation in the
manifold case, because this gives the motivation for the different
approaches in the case of fractals. Let $M$ be a compact Riemannian manifold
with a Riemannian metric $g$ given as a quadratic form $g_x$ on the
tangent space $T_xM$ for $x\in M$. As usual, we assume that the
dependence of $g_x$ on $x$ is differentiable. Then the quadratic form
$g_x$ defines an isomorphism $\alpha_g$ between the tangent space
$T_xM$ and its dual $T^*_xM$ (and thus on the tangent bundle $TM$ and
the cotangent bundle $T^*M$) by $\alpha_g(v)w=g_x(v,w)$ for $v,w\in
T_xM$. This defines the gradient of a function $f$ as $\grad
f=\alpha_g^{-1}(df)$. Define the divergence of a vector field $X$ as
the negative formal adjoint of $\grad$ with respect to the scalar product
$\langle X,Y\rangle_{L^2(M)}=\int_M g_x(X,Y)\,\dd\vol(x)$:
\begin{equation*}
  \langle X,\grad f\rangle_{L^2(M)}=-\langle \divg X,f\rangle_{L^2(M)}.
\end{equation*}
The Laplace operator is then defined as
(cf.~\cite{Rosenberg1997:laplacian_riemannian_manifold,
  Federer1969:geometric_measure_theory})
\begin{equation}\label{eq:divgrad}
\Laplace f=\divg\grad f.
\end{equation}
By definition this operator is self-adjoint and thus has only non-positive real
eigenvalues by $\langle \Laplace f,f\rangle_{L^2}=-\langle\grad f,\grad
f\rangle_{L^2}\leq0$.

Based on the above approach, a corresponding energy form (``Dirichlet form'',
cf.~\cite{Fukushima_Oshima_Takeda2011:dirichlet_forms_markov}) can be defined
\begin{equation*}
\mathcal{E}(u,v)=\int_M g_x(\grad u,\grad v)\,\dd\vol(x)=
\int_M g_x(\dd u,\dd v)\,\dd\vol(x),
\end{equation*}
which lends itself to a further way of defining a Laplace operator via the
relation 
\begin{equation}\label{eq:mathc-v=langle-lapl}
\mathcal{E}(u,v)=-\langle \Laplace u,v\rangle_{L^2(M)}.
\end{equation}

Geometrically, the Laplace operator measures the deviation of the
function $f$ from the mean value. More precisely, let
$S(x,r)=\{y\in M\mid d(x,y)=r\}$ denote the ball of radius $r$ in $M$
(in the Riemannian metric). Then
\begin{equation}\label{eq:lapl-mean}
\Laplace f=2n\lim_{r\to0+}\frac1{r^2\sigma(S(x,r))}\int_{S(x,r)}(f(y)-f(x))\,
\dd\sigma(y),
\end{equation}
where $n$ denotes the dimension of the manifold $M$ and $\sigma$ is
the surface measure on $S(x,r)$. This also motivates the definition of the
Laplace operator as limit of finite difference operators
\begin{equation*}
\Laplace f(x)=\lim_{r\to0}\frac1{r^2}\left(\sum_{p\in N_r(x)}w_pf(p)-f(x)\right),
\end{equation*}
where $N_r(x)$ is a finite set of points at distance $r$ from $x$ and $w_p$ are
suitably chosen weights. Such approximations to the Laplace operator are the
basis of the method of finite differences in numerical mathematics.

The Laplace operator can then be used to define a diffusion on $M$ via the heat
equation $\Laplace u=\partial_t u$. The solution $u(t,x)$ of the initial value
problem $u(0,x)=f(x)$ defines a semi-group of operators
$P_t$ by
\begin{equation*}
u(t,x)=P_tf(x).
\end{equation*}
The semi-group property $P_{s+t}=P_sP_t$ comes from the uniqueness of the
solution $u$ and translation invariance with respect to $t$ of the heat
equation. From the heat semi-group $P_t$ the Laplace operator can be recovered
as the infinitesimal generator
\begin{equation}\label{eq:lapl-generator}
\Laplace f=\lim_{t\to0+}\frac{P_tf-f}t,
\end{equation}
which exists on a dense subspace of $L^2(M)$ under suitable continuity
assumptions on the semi-group $(P_t)_{t\geq0}$
(cf.~\cite{Yosida1971:functional_analysis}).

In the fractal situation none of the above approaches can be used directly to
define a Laplace operator. The main reason for this is that there is no natural
definition of derivative on a fractal. But the above approaches to the Laplace
operator on a manifold can be used in the opposite direction: 
\begin{itemize}
\item starting from a diffusion process that can be defined on fractals by
  approximating random walks. Then the Laplace operator can be defined as the
  infinitesimal generator. This is described in
  Sections~\ref{sec:random-walks-graphs} and~\ref{sec:laplace-operator-inf}.
\item taking the limit of finite difference operators on graphs approximating
  the fractal gives a second possible approach to the Laplace operator, which
  is explained in Section~\ref{sec:laplace-operator-limit}.
\item starting with a Dirichlet form $\mathcal{E}$ gives a third possible
  approach, which is presented in Section~\ref{sec:lapl-oper-dirichlet}.
\end{itemize}
\subsection{Random walks on graphs and diffusion on 
fractals}\label{sec:random-walks-graphs}
The first idea to define a diffusion on a fractal was to define a
sequence of random walks on approximating graphs and to synchronise
time so that the limiting process is non-constant and continuous. This
was the first approach to the diffusion process on the Sierpi\'nski
gasket given in \cite{Goldstein1987:random_walks_diffusion,
  Kusuoka1987:diffusion_process_fractal,
  Barlow_Perkins1988:brownian_motion_sierpinski} and later generalised
to other ``nested fractals'' in
\cite{Lindstroem1990:brownian_motion_nested}. Because of its
importance for our exposition, we will explain it in some detail in
this section. We will follow the lines of definition of self-similar
graphs given in \cite{Kroen2002:green_functions_self,
  Kroen_Teufl2004:asymptotics_transition_probabilities} and adapt it
for our purposes.

We consider a graph $G=(V(G),E(G))$ with vertices $V(G)$ and undirected edges
$E(G)$ denoted by $\{x,y\}$. We assume throughout that $G$ does not contain
multiple edges nor loops. For $C\subset V(G)$ we call $\partial C$ the vertex
boundary, which is given by the set of vertices in $V(G)\setminus C$, which are
adjacent to a vertex in $C$. For $F\subset V(G)$ we define the reduced graph
$G_F$ by $V(G_F)=F$ and $\{x,y\}\in E(G_F)$, if $x$ and $y$ are in the boundary
of the same component of $V(G)\setminus F$.
\begin{defn}\label{def:self-similar}
A connected infinite graph $G$ is called self-similar with respect to $F\subset
V(G)$ and $\phi:V(G)\to V(G_F)$, if
\begin{enumerate}
\item\label{decim1} no vertices in $F$ are adjacent in $G$
\item the intersection of the boundaries of two different components of
  $V(G)\setminus F$ does not contain more than one point
\item $\phi$ is an isomorphism of $G$ and $G_F$.
\end{enumerate}
\end{defn}
A random walk on $G$ is given by transition probabilities $p(x,y)$, which are
positive, if and only if $\{x,y\}\in E(G)$. For a trajectory $(Y_n)_{n\in\N_0}$
of this random walk with $Y_0=x\in F$ we define stopping times recursively by
\begin{equation*}
T_{m+1}=\min\left\{k>T_m\mid Y_k\in F\setminus\{Y_{T_m}\}\right\},\quad T_0=0.
\end{equation*}
Then $(Y_{T_m})_{m\in\N_0}$ is a random walk on $G_F$. Since the underlying
graphs $G$ and $G_F$ are isomorphic, it is natural to require that
$(\phi^{-1}(Y_{T_m}))_{m\in\N_0}$ is the same stochastic process as
$(Y_n)_{n\in\N_0}$. This requires the validity of equations for the basic
transition probabilities
\begin{equation}\label{eq:transition}
\P\left(Y_{T_{n+1}}=\phi(y)\mid Y_{T_m}=\phi(x)\right)=
\P\left(Y_{n+1}=y\mid Y_n=x\right)=p(x,y).
\end{equation}
These are usually non-linear rational equations for the transition
probabilities $p(x,y)$. The existence of solutions of these equations has been
the subject of several investigations, and we refer to
\cite{Metz1995:hilberts_projective_metric, Metz1996:renormalization_contracts,
  Metz2003:cone_diffusions,Sabot1997:existence_uniqueness_diffusion}.

The process $(Y_n)_{n\in\N_0}$ on $G$ and its ``shadow'' $(Y_{T_n})_{n\in\N_0}$
on $G_F$ are equal, but they are on a different time scale. Every transition
$Y_{T_n}\to Y_{T_{n+1}}$ on $G_F$ comes from a path $Y_{T_n}\to
Y_{T_n+1}\cdots\to Y_{T_{n+1}-1}\to Y_{T_{n+1}}$ in a component of $V(G)\setminus
F$. The time scaling factor between these processes is given by
\begin{equation*}
\lambda=\E(T_{n+1}-T_n)=\E(T_1).
\end{equation*}
This factor is $\geq2$ by assumption \eref{decim1} on $F$. More precisely, the
relation between the transition time on $G_F$ and the transition time on $G$ is
given by a super-critical ($\lambda>1$) branching process, which replaces an
edge $\{\phi(x),\phi(y)\}\in G_F$ by a path in $G$ connecting the points $x$
and $y$ without visiting a point in $V(G)\setminus F$ (except for $x$ and for
$y$ in the last
step).

In order to obtain a process on a fractal in $\R^d$, we assume further that $G$
is embedded in $\R^d$ (i.e. $V(G)\subset\R^d$). The self-similarity of the
graph is carried over to the embedding by assuming that there exists a
$\beta>1$ (the space scaling factor) such that $F=V(G_F)=\beta V(G)$. The
fractal limiting structure is then given by
\begin{equation*}
  Y_G=\overline{\bigcup_{n=0}^\infty\beta^{-n}V(G)}.
\end{equation*}

Iterating this graph decimation we obtain a sequence of (isomorphic)
graphs $G_n=(\beta^{-n}V(G),E(G))$ on different scales. The random
walks $(Y_k^{(n)})_{k\in\N_0}$ on $G_n$ are connected by time scales
with the scaling factor $\lambda$. From the theory of branching
processes (cf.~\cite{Harris1963:theory_branching_processes}) it
follows that the time on level $n$ scaled by $\lambda^{-n}$ tends to a
random variable. From this it follows that
$\beta^{-n}Y_{[t\lambda^n]}$ weakly tends to a (continuous time)
stochastic process $(X_t)_{t\geq0}$ on the fractal $Y_G$.

\begin{figure}[h]
  \centering
  \includegraphics[width=0.8\textwidth]{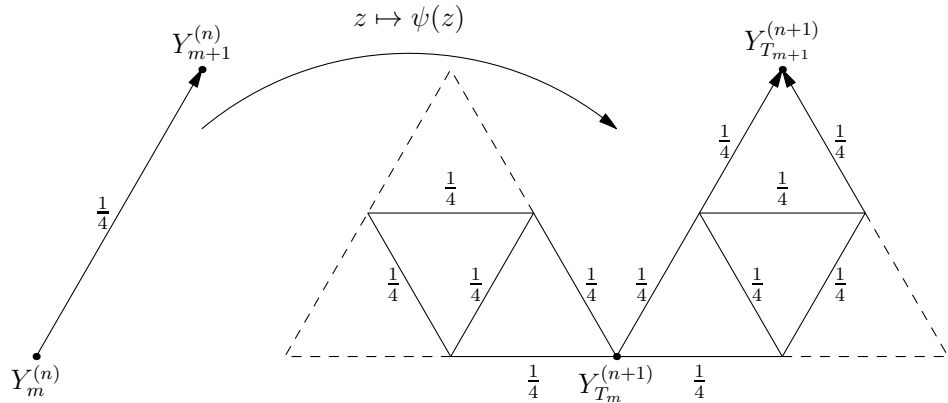}
  \caption{Transition between $Y_k^{(n)}$ and $Y_k^{(n+1)}$}
  \label{fig:psi}
\end{figure}

Under the assumption that the graph $G$ has an automorphism group
acting doubly transitively on the points of the boundary of every
component of $G\setminus F$ (cf.~\cite{Kroen2002:green_functions_self,
  Kroen_Teufl2004:asymptotics_transition_probabilities}), we consider
the \emph{generating function} of the transition probabilities
$p_n(x,y)=\P(Y_n=y\mid Y_0=x)$, the so called \emph{Green function} of
the random walk
\begin{equation*}
G(z\mid x,y)=\sum_{n=0}^\infty p_n(x,y)z^n,
\end{equation*}
which can also be seen as the resolvent $(I-zP)^{-1}$ of the transition
operator $P$. The replacement rules connecting the random walks on the graphs
$G$ and $G_F$ result in a functional equation for the Green function
\begin{equation}\label{eq:green}
G(z\mid\phi(x),\phi(y))=f(z)G(\psi(z)\mid x,y),
\end{equation}
where the rational function $\psi(z)$ is the probability generating function of
the paths connecting two points $v_1,v_2\in\phi(V(G))$ without reaching any
point in $\phi(V(G))\setminus\{v_1\}$.  Then the time scaling factor is the
expected number of steps needed for the paths counted by $\psi$, so we have
\begin{equation*}
  \lambda=\psi'(1).
\end{equation*}
The rational function $f(z)$ is the
probability generating function of the paths starting and ending in
$v_1\in\phi(V(G))$ without reaching any other point in
$\phi(V(G))$. Equation \eref{eq:green} becomes especially simple for a fixed
point $x$ of the map $\phi$
\begin{equation*}
G(z\mid x,x)=f(z)G(\psi(z)\mid x,x).
\end{equation*}
It was proved in \cite{Kroen_Teufl2004:asymptotics_transition_probabilities}
that under our conditions on the set $F$ the map $\phi$ can have at most one
fixed point.

The Koenigs function $\Phi$ of $\psi$ around the fixed point $z=1$ is given by
\begin{equation*}
\Phi(\lambda z)=\psi(\Phi(z)),\quad \Phi(0)=1,\quad \Phi'(0)=1.
\end{equation*}
This can be used to linearise this functional equation and to obtain precise
analytic information about $G(z\mid x,x)$
(cf.~\cite{Kroen_Teufl2004:asymptotics_transition_probabilities,
  Teufl2007:asymptotic_behaviour_analytic}). From this the asymptotic behaviour
of the transition probabilities $p_n(x,x)$ can be derived. In many examples
these transition probabilities exhibit periodic fluctuations
\begin{equation*}
p_n(x,x)\sim n^{-\frac{d_S}2}\left(\sigma(\log_\lambda n)+\BigOh(n^{-1})\right),
\end{equation*}
where $\sigma$ is a continuous, periodic, non-constant function of period $1$
(cf.~\cite{Kroen_Teufl2004:asymptotics_transition_probabilities}).

A first example of such graphs is the Sierpi\'nski graph studied as an
approximation to the fractal Sierpi\'nski gasket. In this case we have for the
probability generating function $\psi(z)=\frac{z^2}{4-3z}$ and $\lambda=5$. The
random walk on this graph was studied in
\cite{Barlow_Perkins1988:brownian_motion_sierpinski} in order to define a
diffusion on a fractal set. Self similarity of the graph and the fractal have
been exploited further, to give a more precise description of the random walk
\cite{Grabner_Woess1997:functional_iterations_periodic} and the diffusion
\cite{Grabner1997:functional_iterations_stopping}. In
\cite{Kroen_Teufl2004:asymptotics_transition_probabilities} a precise
description of the class of graphs in terms of their symmetries is given, which
allow a similar construction. In \cite{Schweitzer2006:diffusion_simply_doubly}
this analysis was carried further to obtain results for the transition
probabilities under less symmetry assumptions using multivariate generating
functions.

In some examples (for instance the Sierpi\'nski graph) it occurred that the
function $\psi$ was conjugate to a polynomial $p$, i.e.
\begin{equation*}
\psi(z)=\frac1{p(\frac1z)},
\end{equation*}
which allowed a study of the properties of the random walk by
referring to the classical Poincar\'e equation. Properties of the
Poincar\'e equation have been used in
\cite{Derfel_Grabner_Vogl2008:zeta_function_laplacian} to study the
analytic properties of the zeta function of the Laplace operator given
by the diffusion on certain self-similar fractals. For more details we
refer to Section~\ref{sec:spectr-zeta-funct}.

\begin{remark}
  The Sierpi\'nski carpet is a typical example of an infinitely
  ramified fractal. In
  \cite{Barlow_Bass1999:graphical_sierpinski_carpet} approximations by
  graphs are used to define a diffusion on this fractal. By the
  infinite ramification, this approach is is more intricate than the
  procedure described here. By the very recent result on the
  uniqueness of Brownian motion on the Sierpi\'nski carpet
  (cf.~\cite{Barlow_Bass_Kumagai+2010:uniqueness_of_brownian}) this
  yields the same process as constructed by rescaling the classical
  Brownian motion restricted to finite approximations of the
  Sierpi\'nski carpet in
  \cite{Barlow_Bass1989:construction_brownian_motion,
    Barlow_Bass1999:brownian_sierpinski_carpet}
\end{remark}
\begin{figure}[h]
  \centering
  \includegraphics[width=0.5\textwidth]{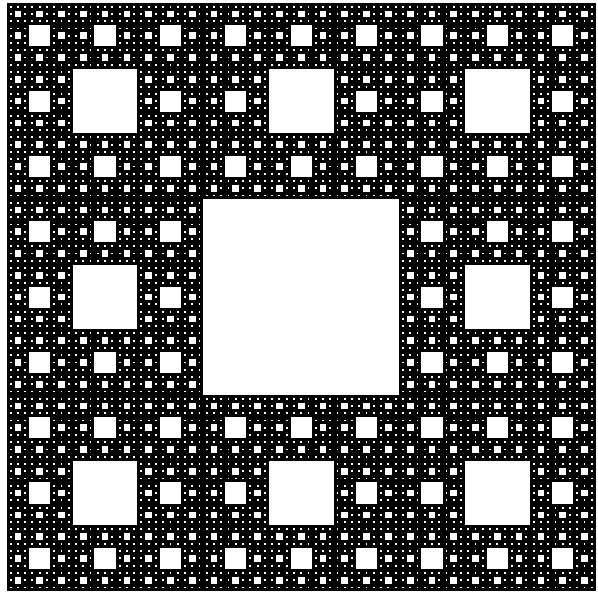}
  \caption{The Sierpi\'nski carpet}
  \label{fig:carpet}
\end{figure}

\subsection{The Laplace operator as the infinitesimal 
generator of a diffusion}\label{sec:laplace-operator-inf}
Given a diffusion process $(X_t)_{t\geq0}$ on a fractal $K$ we can now define a
corresponding Laplace operator. At first we define a semi-group of operators
$A_t$ by
\begin{equation*}
A_tf(x)=\E_x f(X_t)
\end{equation*}
for functions $f\in L^2(K)$. The semi-group property
\begin{equation*}
A_sA_t=A_{s+t}
\end{equation*}
of the operators comes from the Markov property of the underlying stochastic
process $X_t$.

By \cite[Chapter~IX]{Yosida1971:functional_analysis} this semi-group has an
infinitesimal generator given as
\begin{equation*}
\Delta f=\lim_{t\to0+}\frac{A_tf-f}t.
\end{equation*}
This limit exists on a dense subspace $\F$ of $L^2(K)$ and is called Laplace
operator on $K$. This name comes from the fact that for the usual Brownian
motion on a manifold this procedure yields the classical Laplace-Beltrami
operator. The function $u(x,t)=A_tf(x)$ satisfies the heat equation
\begin{equation*}
\Laplace u=\partial_tu,\quad u(x,0)=f(x).
\end{equation*}

It was observed in the early beginnings of the development of the theory of
diffusion of fractals that the domain of $\Delta$ does not contain the
restriction of any non-constant differentiable function
(cf.~\cite{Ben-Bassat_Strichartz_Teplyaev1999:not_domain_laplacian}).

\subsection{The Laplace operator as limit of difference
operators on graphs}\label{sec:laplace-operator-limit}
A totally different and more direct approach to the Laplace operator on
self-similar fractals has been given by Kigami in
\cite{Kigami1989:harmonic_calculus_sierpinski}. The operator $\Delta$ is
approximated by difference operators on the approximating graphs $G_n$. The
graph Laplacians are given by
\begin{equation*}
\Delta_n f(x)=\sum_{y\sim x, y\in G_n}p(x,y)f(y)-f(x)
\end{equation*}
as a weighted sum over the neighbours of $x$ in $G_n$ ($y\sim x$ describes the
neighbourhood relation in the graph $G_n$). In order to make this
construction compliant with the approach via stochastic processes, these
operators have to be rescaled appropriately. The correct rescaling is then
given by the time scaling factor $\lambda$ introduced before, namely
\begin{equation*}
\Delta f(x)=\lim_{n\to\infty}\lambda^n\Delta_n f(x).
\end{equation*}

\subsection{Laplace operators via Dirichlet forms}
\label{sec:lapl-oper-dirichlet}
Following the exposition in 
\cite[Chapters~2 and~3]{Kigami2001:analysis_fractals} we define a sequence of 
quadratic forms on the finite sets $V_m$ given in Definition~\ref{def:vm}.
\begin{defn}\label{def:dirichlet}
  Let $V$ be a finite set. Then a bilinear form $\mathcal{E}$ on $\ell(V)$, the
  real functions on $V$ is called a Dirichlet form, if the following conditions
  hold:
\begin{enumerate}
\item $\forall u\in\ell(V):\mathcal{E}(u,u)\geq0$
\item $\mathcal{E}(u,u)=0$ implies that $u$ is constant on $V$
\item for $u\in\ell(V)$ and $\overline{u}(x)=\max(0,\min(u(x),1))$ the
  inequality $\mathcal{E}(u,u)\geq\mathcal{E}(\overline{u},\overline{u})$
  holds.
\end{enumerate}
\end{defn}

\begin{defn}
Let $\mathcal{E}$ be a Dirichlet form on the finite set $V$ and let $U$ be a
proper subset of $V$. Then the restriction of $\mathcal{E}$ to $U$ is defined
as 
\begin{equation}\label{eq:restr}
\mathcal{R}_{V,U}(\mathcal{E})(u,u)=
\min\left\{\mathcal{E}(v,v)\mid v\in\ell(V), v|_U=u\right\}.
\end{equation}
\end{defn}
On the level of the coefficient matrices of the Dirichlet forms, the
operation of restriction is given in terms of the Schur complement.
\begin{defn}\label{def:compatible}
Let $(V_n,\mathcal{E}_n)_n$ be a sequence of increasing finite sets $V_n$ and
Dirichlet forms $\mathcal{E}_n$ on $V_n$. The sequence is called
\emph{compatible}, if
\begin{equation*}
\mathcal{R}_{V_{n+1},V_n}(\mathcal{E}_{n+1})=\mathcal{E}_n
\end{equation*}
holds for all $n$.
\end{defn}

For a compatible sequence $(V_n,\mathcal{E}_n)_n$ and a function $f$ on $K$,
the sequence $(\mathcal{E}_n(f|_{V_n},f|_{V_n}))_n$ is increasing by definition
and thus converges to a value in $[0,\infty]$. This makes the following
definition natural.
\begin{defn}\label{def:E}
Let $(V_n,\mathcal{E}_n)_n$ be a compatible sequence of Dirichlet forms. Let
\begin{equation*}
\mathcal{D}=\left\{f:K\to\R\mid 
\lim_{n\to\infty}\mathcal{E}_n(f|_{V_n},f|_{V_n})<\infty\right\}.
\end{equation*}
Then for all $f\in\mathcal{D}$
\begin{equation*}
\mathcal{E}(f,f)=\lim_{n\to\infty}\mathcal{E}_n(f|_{V_n},f|_{V_n})
\end{equation*}
defines a Dirichlet form on $K$, and $\mathcal{D}$ is its domain.
\end{defn}

In order to make the sequence of Dirichlet forms coherent with the self-similar
structure of $K$, we require the following self-similarity condition for
$\mathcal{E}_n$
\begin{equation*}
\mathcal{E}_{n+1}(f,f)=
\lambda\sum_{i=1}^m r_i^{-1} \mathcal{E}_n(f\circ F_i,f\circ F_i),
\end{equation*}
where $r_i$ ($i=1,\ldots,m$) are positive weights and $\lambda$ is a
proportionality factor. Furthermore, the sequence of forms has to be
compatible, which amounts to the equation
\begin{equation}\label{eq:compat0}
\lambda\mathcal{R}_{V_1,V_0}\left(\sum_{i=1}^mr_i^{-1}
\mathcal{E}_0(\cdot\circ F_i,\cdot\circ F_i))\right)=\mathcal{E}_0(\cdot,\cdot),
\end{equation}
which comes as a solution of a non-linear eigenvalue equation. This equation
plays the same role in the Dirichlet form approach as the equations
\eref{eq:transition} play for the approach via random walks. The Dirichlet form
$\mathcal{E}$ on $K$ is then defined as
\begin{equation*}
\mathcal{E}(f,f)=\lim_{n\to\infty}\lambda^{n}\sum_{w\in S^n}r_w^{-1}
\mathcal{E}_0(f\circ F_w|_{V_0},f\circ F_w|_{V_0}),
\end{equation*}
where $S=\{1,\ldots,m\}$, $r_w=r_{w_1}\cdots r_{w_n}$ for $w=w_1\ldots w_n$ and
$F_w=F_{w_1}\circ \cdots\circ F_{w_n}$.

\begin{remark}
There are some additional technical problems concerning this construction of
Dirichlet forms, which arise from the fact that in general the form is
supported only on a proper subset of $K$. In
\cite[Chapter~3]{Kigami2001:analysis_fractals} sufficient conditions for the
weights $r_i$ and the form $\mathcal{E}_0$ are given, which ensure that the
form $\mathcal{E}$ is supported on the whole set $K$.
\end{remark}

Given a Dirichlet form $\mathcal{E}$ together with its domain $\mathcal{D}$ and
a measure $\mu$ on $K$, we can now define the associated Laplace operator on
$K$ by
\begin{equation}\label{eq:laplace}
\forall v\in\mathcal{D}\cap L^2(\mu):
\mathcal{E}(u,v)=-\langle {\Laplace}_\mu u,v\rangle_{L^2(\mu)},
\end{equation}
which defines $\Laplace_\mu$, if $\mathcal{D}$ is dense in $L^2(\mu)$. Notice
that this is the same equation as in the manifold case
\eref{eq:mathc-v=langle-lapl}.

In the case of a self-similar fractal $K$ as described in
Section~\ref{sec:fractals} the ``natural'' measure on $K$ is the according
Hausdorff measure $\mathcal{H}_K^{\rho}$. In this case we omit the subscript
$\mu$. 

For more information on Dirichlet forms in general and their applications to
the description of diffusion processes, we refer to the monograph
\cite{Fukushima_Oshima_Takeda2011:dirichlet_forms_markov}. For the specifics of
Dirichlet forms on fractals we refer to 
\cite{Kigami2001:analysis_fractals,Fukushima1992:dirichlet_forms}.



%% file: Spectral.tex
\section{Spectral analysis on fractals}
\label{sec:spectr-analys-fract}
\subsection{Spectral analysis on manifolds}
\label{sec:spectral_manifolds}
Let us start with a short discussion of the manifold case, which is again
somehow complementary to the fractal case. As described in
Section~\ref{sec:laplace} the Laplace operator on a manifold $M$ is defined via
a Riemannian metric $g$. If $M$ is compact with smooth (or empty) boundary (for
simplicity) the Laplace operator has pure point spectrum. The eigenvalues
$-\lambda_k$ are all real (by self-adjointness) and non-positive (by
definition \eref{eq:divgrad}). Denote the normalised eigenfunctions by
$\psi_k$. Then the heat kernel can be written as
\begin{equation}\label{eq:heat-kernel}
p_t(x,y)=\sum_{k=0}^\infty e^{-\lambda_k t}\psi_k(x)\psi_k(y).
\end{equation}
This expression yields the trace of the heat kernel
\begin{equation}\label{eq:heat-trace}
K(t)=\int_M p_t(x,x)\,\dd\vol(x)=\sum_{k=0}^\infty e^{-\lambda_k t}.
\end{equation}

On the other hand, if $M$ is closed, the asymptotic behaviour of the heat
kernel for $t\to0+$ can be described very precisely
\begin{equation}\label{eq:heat-asymptotic}
p_t(x,y)=(4\pi t)^{-\frac n2}e^{-d(x,y)^2/(4t)}
\sum_{\ell=0}^\infty t^\ell a_\ell(x,y),
\end{equation}
where the functions $a_\ell$ can be computed iteratively by solving certain
second order partial differential equations
(cf.~\cite{Minakshisundaram_Pleijel1949:some_properties_eigenfunctions,
Berline_Getzler_Vergne1992:heat_kernels}). In particular,
$a_0(x,y)=1$. The values $a_\ell(x,x)$ can be expressed in terms of the
functions $\left.\Laplace_y^kd(x,y)^{2m}\right|_{y=x}$
(cf.~\cite{Polterovich2000:heat_invariants}). In the case of manifolds $M$ with
smooth boundary, further terms involving half integer powers of $t$ occur
\begin{equation}\label{eq:heat-boundary}
p_t(x,x)=t^{-\frac n2}\sum_{\ell=0}^\infty t^\ell a_\ell(x,x)+
t^{-\frac n2}\sum_{\ell=1}^\infty t^{\frac\ell2}b_\ell(x,x);
\end{equation}
the terms $b_\ell$ encode curvature information of the boundary $\partial
M$. For more details we refer to
\cite{Kirsten2002:spectral_functions_mathematics}.

Especially, this gives an asymptotic expansion of $K(t)$ (in the case of a
closed manifold)
\begin{equation}\label{eq:trace-asymptotic}
K(t)=(4\pi t)^{-\frac n2}\sum_{\ell=0}^\infty t^\ell 
\int_M a_\ell(x,x)\,\dd\vol(x).
\end{equation}
From the first order asymptotic relation $K(t)\sim \vol(M)/(4\pi t)^{\frac n2}$
the asymptotic behaviour of the counting function
\begin{equation}\label{eq:counting}
N(x)=\sum_{\lambda_k<x}1
\end{equation}
can be obtained by a Tauberian argument giving the classical Weyl asymptotic
relation
\begin{equation}\label{eq:weyl-classic}
N(x)=\frac{\vol(M)}{B_n}x^{\frac n2}+o(x^{\frac n2}),
\end{equation}
where $B_n$ denotes the volume of the $n$-dimensional unit ball. The
eigenvalues of $\Laplace$ constitute the frequencies of the oscillations in the
solutions of the wave equation $\Laplace u=u_{tt}$. This led to M.~Kac's
famous question ``Can one hear the shape of a drum?''
(cf.~\cite{Kac1966:hear_drum}).

The precise order of magnitude of the error term has been determined as
$x^{\frac{n-1}2}$ in the case of a smooth boundary of $M$ in
\cite{Ivrii1980:laplace_beltrami}. In the case of fractal boundary, this is
subject of the Weyl-Berry conjecture stated in
\cite{Berry1980:wave_motion,Berry1979:distribution_modes_fractal}. The original
conjecture estimated the error term by $\mathcal{O}(x^{\frac {d_H(\partial
    M)}2})$, where $d_H(\partial M)$ denotes the Hausdorff dimension of the
boundary of $M$. This was shown to be false in
\cite{Brossard_Carmona1986:can_one_hear}. The Hausdorff dimension was then
replaced by the Minkowski dimension $d_M(\partial M)$. It was shown in
\cite{Lapidus1991:fractal_drum_inverse} that the error term is
$\mathcal{O}(x^{\frac d2})$ for any $d>d_M(\partial M)$. Furthermore, it was
proved in \cite{Lapidus_Pomerance1993:riemann_weyl_berry} that the modified
conjecture is true for dimension $n=1$; in
\cite{Lapidus_Pomerance1996:counterexamples_weyl_berry} counterexamples for
dimension $n>1$ are presented. In these counterexamples disconnected manifolds
are considered, thus the conjecture could still hold for connected manifolds
$M$ with fractal boundary.

More precise information on the eigenvalues is contained in the spectral
zeta function
\begin{equation}\label{eq:zeta}
\zeta_{\Laplace}(s)=\sum_{\lambda_k\neq0}(\lambda_k)^{-s},
\end{equation}
the Dirichlet generating function of $(\lambda_k)_k$. The zeta function is
connected to the trace of the heat-kernel by a Mellin transform
\begin{equation}\label{eq:Mellin-trace}
\zeta_{\Laplace}(s)\Gamma(s)=\int_0^\infty K(t)t^{s-1}\,\dd t, \quad\mbox{ for }
\Re s>\frac n2.
\end{equation}
The right hand side has an analytic continuation to the whole complex plane,
which can be found by using the asymptotic expansion
\eref{eq:trace-asymptotic}. Furthermore, (in the case of a closed manifold) the
right hand side has at worst simple poles at the points $s=\frac n2-\ell$ (with
$\ell\in\{0,1,\ldots, \lfloor\frac{n-1}2\rfloor$ for even $n$ and $\ell\in\N_0$
for odd $n$) with residues
\begin{equation*}
\Res_{s=\frac n2-\ell}\zeta_{\Laplace}(s)=
\frac1{\Gamma(\frac n2-\ell)}\int_M a_\ell(x,x)\,\dd\vol(x)
\end{equation*}
and special values (for $\ell\in\N$)
\begin{equation*}
\zeta_{\Laplace}(-\ell)=
\left\{\begin{array}{ll}
0&\mbox{for }n\mbox{ odd}\\
(-1)^\ell\ell!\int_M a_\ell(x,x)\,\dd\vol(x)&\mbox{for }n\mbox{ even}
\end{array}\right.
\end{equation*}
and
\begin{equation*}
\zeta_{\Laplace}(0)=\left\{
\begin{array}{ll}
-\dim\ker\Laplace&\mbox{for }n\mbox{ odd}\\
\int_M a_{n/2}(x,x)\,\dd\vol(x)-\dim\ker\Laplace&\mbox{for }n\mbox{ even}\\
\end{array}\right.
\end{equation*}
(cf.~\cite{Minakshisundaram_Pleijel1949:some_properties_eigenfunctions}).  The
special value $\zeta'_{\Laplace}(0)$ can be interpreted as the negative
logarithm of the determinant of $\Laplace$.  In the case of a manifold $M$ with
smooth boundary a similar reasoning gives simple poles of the zeta function at
the negative half-integers $-\frac12,-\frac32,\ldots$ and at the points
$\frac n2-\frac12,\frac n2-\frac32,\ldots, \frac{3-(-1)^n}4$ with residues
depending on the functions $b_\ell$.

In the case of a fractal, the above approach to the spectral zeta function, its
analytic continuation, and its finer properties can not be used. The reason for
this is that no asymptotic expansion of the heat kernel is known in
the fractal case. In Section~\ref{sec:trace-heat-kernel} we comment on the
known upper and lower estimates for the heat kernel.

On the other hand, for fractals having spectral decimation, the eigenvalues can
be described very precisely. In this case they turn out to constitute a finite
union of level sets of solutions of the classical Poincar\'e functional
equation. This allows then to obtain the analytic continuation of the spectral
zeta function from the precise knowledge of the asymptotic behaviour of the
Poincar\'e function; the asymptotic expansion of the trace of the heat kernel
can then be obtained from the poles of the zeta function, reversing the
argument in \eref{eq:Mellin-trace}. The eigenvalue counting function can then
be related to the harmonic measure on the Julia set of a the polynomial
governing the spectral decimation.
\subsection{Spectral decimation}\label{sec:spectral-decimation}
It has been first observed by Fukushima and Shima
\cite{Fukushima_Shima1992:spectral_analysis_sierpinski,
  Shima1993:eigenvalue_problem_laplacian,
  Shima1996:eigenvalue_problems_laplacians} that the eigenvalues of
the Laplacian on the Sierpi\'nski gasket and its higher dimensional
analogues exhibit the phenomenon of \emph{spectral
  decimation}. Lateron, spectral decimation for more general fractals
has been studied by Malozemov, Strichartz, and Teplyaev
\cite{Malozemov_Teplyaev2003:self_similarity_operators,
  Strichartz2003:fractafolds_spectra,
  Teplyaev2007:spectral_zeta_functions}.

\begin{defn}[Spectral decimation]\label{def1}
  The Laplace operator on a p.~c.~f.\ self-similar fractal $G$ admits
  \emph{spectral decimation}, if there exists a rational function $R$, a finite
  set $A$ and a constant $\lambda>1$ such that all eigenvalues of $\Delta$ can
  be written in the form 
\begin{equation}\label{eq:decimation}
\lambda^m\lim_{n\to\infty}\lambda^nR^{(-n)}(\{w\}),\quad w\in A, m\in\N
\end{equation}
where the preimages of $w$ under $n$-fold iteration of $R$ have to be chosen
such that the limit exists. Furthermore, the multiplicities $\beta_m(w)$ of the
eigenvalues depend only on $w$ and $m$, and the generating functions of the
multiplicities are rational.
\end{defn}
The fact that all eigenvalues
of $\Delta$ are negative real implies that the Julia set of $R$ has to be
contained in the negative real axis. We will exploit this fact later.

In many cases such as the higher dimensional Sierpi\'nski gaskets, the
rational function $R$ is conjugate to a polynomial. The method for
meromorphic continuation of $\zeta_\Delta$ given in
Section~\ref{sec:spectr-zeta-funct} makes use of this assumption.
Recentlty, Teplyaev \cite{Teplyaev2007:spectral_zeta_functions} showed
under the same assumption that the zeta function of the Laplacian
admits a meromorphic continuation to $\Re s>-\eps$ for some $\eps>0$
depending on properties of the Julia set of the polynomial given by
spectral decimation. His method uses ideas similar to those used in
\cite{Grabner1997:functional_iterations_stopping} for the meromorphic
continuation of a Dirichlet series attached to a
polynomial. Complementary to the ideas used here, Teplyaev's method
carries over to rational functions $R$.

\subsection{Eigenvalue counting}\label{sec:eigenvalue-counting}
As in the euclidean case the eigenvalue counting function
\begin{equation*}
N(x)=\sum_{\lambda_k<x}1
\end{equation*}
measures the number of eigenvalues less than $x$. In a physical context this
quantity is referred to as the ``integrated density of states''. It turns out
that $N(x)$ in many cases does not exhibit a pure power law as in the euclidean
case, but shows periodic fluctuations. One source of this periodicity
phenomenon is actually spectral decimation, especially the high multiplicities
of eigenvalues, as will become clear in Section~\ref{sec:spectr-zeta-funct}.

Recall the definition of harmonic measure on the Julia set of a polynomial $p$
of degree $d$
(cf.~\cite{Ransford1995:potential_theory_complex_plane}): the sequence of
measures
\begin{equation*}
\mu_n=\frac1{d^n}\sum_{p^{(n)}(x)=\xi}\delta_x
\end{equation*}
converges weakly to a limiting measure $\mu$, the harmonic measure on the Julia
set of $p$. The point $\xi$ can be chosen arbitrarily.

Assume now that the Laplace operator on the fractal $K$ admits spectral
decimation with a polynomial $p$. Then the relation
$\lim_{n\to\infty}\lambda^np^{(-n)}(\{w\})\in B(0,x)$ can be translated into
\begin{equation*}
\lim_{n\to\infty}p^{(n)}(\lambda^{-n}z)=w \mbox{ and } |z|<x.
\end{equation*}
Here $B(0,x)$ denotes the ball of radius $x$ around $0$. By general facts about
polynomial iteration (cf.~\cite{Milnor2006:dynamics_complex}), the limit exists
and defines an entire function of $z$, the Poincar\'e function $\Phi(z)$. The
number of eigenvalues with $m=0$ in $B(0,x)$ is the equal to
\begin{equation*}
N(x)=\sum_{w\in A}N_w(x)
\end{equation*}
with
\begin{equation*}
N_w(x)=\#\left\{z\in B(0,x)\mid \Phi(z)=w\right\}.
\end{equation*}
The following relation can be obtained from the definition of harmonic measure
\begin{equation*}
\lim_{n\to\infty}d^{-n}N_w(\lambda^nx)=\mu\left(\Phi^{-1}(B(0,x))\right);
\end{equation*}
this holds for all $x$ small enough to ensure the existence of the inverse
function $\Phi^{-1}$ on $B(0,x)$.

In \cite[Theorem~5.2]{Derfel_Grabner_Vogl2008:complex_asymptotics_julia} we
could prove a relation between the asymptotic behaviour of the partial counting
functions $N_w(x)$ and the harmonic measure of balls $\mu(B(0,x))$. The
existence of the two limits $(\rho=\log_\lambda d$)
\begin{eqnarray*}
&\lim_{x\to\infty}x^{-\rho}N_w(x)\\
&\lim_{t\to0}t^{-\rho}\mu(B(0,t))
\end{eqnarray*}
is equivalent. We conjectured there, that these limits can only exist, if $p$
is either a Chebyshev polynomial or a monomial. These are the only cases of
polynomials with smooth Julia sets
(cf.~\cite{Hamilton1995:length_julia_curves}).

Summing up the above discussion, the eigenvalue counting function can be
written as
\begin{equation}\label{eq:Nx=sum}
N(x)=\sum_{w\in A}\sum_{m=0}^\infty \beta_m(w)N_w(\lambda^{-m}x).
\end{equation}
Notice, that for fixed $x$, these sums are actually finite. In the known cases,
such as the Sierpi\'nski gasket, the growth of $\beta_m(w)$ is stronger than
$d^m$, which implies that the terms for large $m$ (with $N_w(\lambda^{-m}x)$
still positive) become dominant in this sum. This shows that multiplicity of
the eigenvalues has the main influence on the asymptotic behaviour of
$N(x)$. Furthermore, this explains the presence of an oscillating factor in the
asymptotic main term of $N(x)$. We will discuss that in more detail in
Section~\ref{sec:spectr-zeta-funct}.

%

\subsection{Spectral zeta functions}\label{sec:spectr-zeta-funct}
As in the euclidean case, the eigenvalues of the Laplace operator $\Laplace$
can be put into a Dirichlet generating function. This will later allow to use
methods and ideas from analytic number theory to obtain more precise asymptotic
information on $N(x)$. The zeta function is again given by
\begin{equation*}
\zeta_{\Laplace}(s)=\sum_{\lambda_k\neq0}\lambda_k^{-s},
\end{equation*}
where all eigenvalues are counted with their multiplicity. The zeta function is
related to the eigenvalue counting function by
\begin{equation*}
\zeta_{\Laplace}(s)=\int_0^\infty x^{-s}\dd N(x)=s\int_0^\infty N(x)x^{-s-1}\dd x.
\end{equation*}
The second relation identifies $\zeta_{\Laplace}$ as the Mellin transform of
the counting function $N(x)$.

In the sequel we will exploit the consequences of spectral decimation. Not too
surprisingly after Definition~\ref{def1}, iteration of polynomials will play an
important role in this discussion. Furthermore, since the relation
\eref{eq:decimation} can be expressed in terms of the Poincar\'e function $f$,
properties of this function will be used to derive the meromorphic continuation
of $\zeta_{\Laplace}$ to the whole complex plane.

Under the assumptions of spectral decimation, the Julia set of the polynomial
$p$ is a subset of the non-positive reals, which contains $0$. By
\cite[Theorem~4.1]{Derfel_Grabner_Vogl2008:complex_asymptotics_julia} this
implies that $\lambda=p'(0)\leq d^2$. By
\cite[Theorem~4.1]{Derfel_Grabner_Vogl2008:complex_asymptotics_julia} equality
can only occur, if $p$ is a Chebyshev polynomial, which would correspond to
spectral decimation on the unit interval (viewed as a self-similar
fractal). Thus in the cases of interest, we have that $\rho=\log_\lambda
d<\frac12$. The Poincar\'e function $\Phi$ is then an entire function of order
$\rho$.

In order to find the analytic continuation of $\zeta_{\Laplace}(s)$ to the
whole complex plane, we analyse the partial zeta functions
\begin{equation}\label{eq:partial-zeta}
\zeta_{\Phi,w}(s)=\sum_{\Phi(-\mu)=w,\mu\neq0}\mu^{-s}.
\end{equation}
Since $\Phi_w=1-\frac1w\Phi$ is a function of order $\rho=\log_\lambda
d<\frac12$, it can be expressed as a Hadamard product
\begin{equation*}
1-\frac1w\Phi(z)=
\prod_{\Phi(-\mu)=w}\left(1+\frac z\mu\right);
\end{equation*}
for $w=0$ we have the slightly modified expression
\begin{equation*}
\Phi_0(z)=\frac1z\Phi(z)=
\prod_{\Phi(-\mu)=0,\mu\neq0}\left(1+\frac z\mu\right).
\end{equation*}

Taking the Mellin-transform of $\log\Phi_w$, which exists for $-1<\Re s<-\rho$,
we get
\begin{equation}\label{eq:Mellin-Phiw}
M_w(s)=\int_0^\infty \left(\log\Phi_w(x)\right)x^{s-1}\,\dd x=
\frac\pi{s\sin\pi s}\zeta_{\Phi,w}(-s).
\end{equation}
Thus for finding the analytic continuation of $\zeta_{\Phi,w}(s)$ to the left
of its abscissa of convergence, it suffices to find the analytic continuation
of $M_w(s)$ for $\Re s>-\rho$. Following slightly different lines as in
\cite{Derfel_Grabner_Vogl2008:zeta_function_laplacian}, we consider the
function
\begin{equation*}
\Psi_w(z)=\frac{p(\Phi(z))-w}{a_d(\Phi(z)-w)^d}=
\frac{\Phi_w(\lambda z)}{a_d(-w)^{d-1}\Phi_w(z)^d}.
\end{equation*}
Then we have
\begin{equation*}
\log\Psi_w(z)=\log\Phi_w(\lambda z)-d\log\Phi_w(z)-\log a_d-(d-1)\log(-w)
\end{equation*}
and this function tends to $0$ exponentially for $z\to+\infty$.  Taking the
Mellin transform we obtain
\begin{eqnarray*}
&\fl(\lambda^s-d)M_w(s)\\
&=\int_0^\infty \left(\log\Phi_w(\lambda x)-d\log\Phi_w(x)\right)x^{s-1}\,\dd x\\
&\quad\mbox{for }-1<\Re s<-\rho\\
&=\int_0^1\left(\log\Phi_w(\lambda x)-d\log\Phi_w(x)-\log a_d-
(d-1)\log(-w)\right)x^{s-1}\,\dd x\\
&+(\log a_d+(d-1)\log(-w))\frac1s+
\int_1^\infty\left(\log\Phi_w(\lambda x)-d\log\Phi_w(x)\right)x^{s-1}\,\dd x\\
&\quad\mbox{for }\Re s>-1\\
&=\int_0^1\left(\log\Psi_w(x)\right)x^{s-1}\,\dd x\\
&+\int_1^\infty\left(\log\Phi_w(\lambda x)-d\log\Phi_w(x)-
\log a_d-(d-1)\log(-w)\right)x^{s-1}\,\dd x\\
&\quad\mbox{for }\Re s>0\\
&=\int_0^\infty\left(\log\Psi_w(x)\right)x^{s-1}\,\dd x.
\end{eqnarray*}
Reading the fourth line of this computation shows that $M_w(s)$ has a simple
pole at $s=0$ with residue
\begin{equation*}
\Res_{s=0}M_w(s)=\frac{\log a_d}{d-1}+\log(-w).
\end{equation*}
Furthermore, this computation shows that $M_w(s)$ is holomorphic in the
half-plane $\Re s>0$. Using \eref{eq:Mellin-Phiw} gives the analytic
continuation of $\zeta_{\Phi,w}(s)$ for $\Re s<0$
\begin{equation}\label{eq:cont-zeta}
\zeta_{\Phi,w}(s)=\frac{\lambda^s s\sin\pi s}{\pi(1-d\lambda^s)}
\int_0^\infty\left(\log\Psi_w(x)\right)x^{-s-1}\,\dd x.
\end{equation}
This shows that $\zeta_{\Phi,w}(-m)=0$ for $m\in\N_0$ (for $s=0$ the double
zero of $s\sin\pi s$ cancels the simple pole of $M_w(-s)$). These could be
called the ``trivial zeros'' as in the case of the Riemann zeta
function. Furthermore, we obtain
\begin{equation*}
\zeta'_{\Phi,w}(0)=\frac{\log a_d}{d-1}+\log(-w).
\end{equation*}
The equation \eref{eq:cont-zeta} even lends itself to the numerical
computation of values of $\zeta_{\Phi,w}(s)$ for $\Re s<0$, as we will se in
Section~\ref{sec:numerical}.

By our assumption on spectral decimation, the generating functions of the
multiplicities of the eigenvalues are rational
\begin{equation*}
R_w(z)=\sum_{m=0}^\infty\beta_m(w)z^m.
\end{equation*}
Thus we can write the spectral zeta function of $\Laplace$ as
\begin{equation}\label{eq:zeta=sum}
\zeta_{\Laplace}(s)=
\sum_{w\in A}\sum_{m=0}^\infty\beta_m(w)\sum_{\Phi(-\mu)=w}(\lambda^m\mu)^{-s}=
\sum_{w\in A}R_w(\lambda^{-s})\zeta_{\Phi,w}(s).
\end{equation}
Since all functions involved in the last (finite) sum are meromorphic in the
whole complex plane, we have found the meromorphic continuation of
$\zeta_{\Laplace}$ to the whole complex plane. The functions
$\zeta_{\Phi,w}(s)$ have only simple poles in the points $s=\log_\lambda
d+\frac{2k\pi i}{\log\lambda}$ ($k\in\mathbb{Z}$). Furthermore, the residues of
these poles do not depend on $w$. All other poles of $\zeta_{\Laplace}$ come
from the poles of the functions $R_w(\lambda^{-s})$. Since these are rational
functions of $\lambda^{-s}$, their poles are equally spaced on vertical lines.

Summing up, $\zeta_{\Laplace}(s)$ has a meromorphic continuation to the whole
complex plane with poles in the points
\begin{equation*}
s=-\log_\lambda\beta_{w,j}+\frac{2k\pi i}{\log\lambda}
\mbox{ with }k\in\mathbb{Z},
\end{equation*}
where the $\beta_{w,j}$ are the poles of the rational functions $R_w(z)$ and
(at most) simple poles in the points
\begin{equation*}
s=\log_\lambda d+\frac{2k\pi i}{\log\lambda}
\mbox{ with }k\in\mathbb{Z}.
\end{equation*}
Furthermore, since the functions $\zeta_{\Phi,w}(s)$ are bounded for
$\Re s\geq\log_\lambda d+\eps$ and the functions $R_w(\lambda^{-s})$ are
bounded along every vertical line, which contains no poles, the function
$\zeta_{\Laplace}(s)$ is bounded along every vertical line $c+it$ for
$c>\log_\lambda d$, which does not contain a pole of any of the functions
$R_w(\lambda^{-s})$.

In the case of the Sierpi\'nski gasket and its higher dimensional analogues,
the rightmost poles of $\zeta_{\Laplace}$ come from the rational functions
$R_w(\lambda^{-s})$. Furthermore, the relation
\begin{equation*}
\sum_{w\in A}R_w(1/d)=0
\end{equation*}
holds, which amounts in a (still mysterious) cancellation of the poles of the
functions $\zeta_{\Phi,w}$ in \eref{eq:zeta=sum}. Thus, the analytic behaviour
of the function $\zeta_{\Laplace}$ is mainly governed by the functions
$R_w$. This means that in this respect the effects of multiplicity prevail over
the individual eigenvalues.

We now investigate the asymptotic behaviour of $N(x)$ under the assumption that
$\beta_m(w)$ grows exponentially faster than $d^m$ for some $w\in A$. This
implies that the corresponding term $R_s(\lambda^{-s})$ has poles to the right
of $\Re s=\log_\lambda d$. We use the classical Mellin-Perron formula
(cf.~\cite{Tenenbaum1995:analytic_number_theory}) to
express $N(x)$ in terms of $\zeta_{\Laplace}(s)$
\begin{equation*}
N(x)=\lim_{T\to\infty}\frac1{2\pi i}\int_{c-iT}^{c+iT}
\zeta_{\Laplace}(s)x^s\frac{\dd s}s,
\end{equation*}
for any $c$ such that $R_w(\lambda^{-s})$ has no poles in the half-plane
$\Re(s)\geq c$. Now the line of integration can be shifted to the left to
$\Re(s)=c'$ for $\lambda^{c'}>d$ but such that at least one of the functions
$R_w(\lambda^{-s})$ has poles to the right of $c'$. This is justified, because
$\zeta_{\Laplace}(\sigma+it)$ remains bounded for $|t|\to\infty$ and
$\sigma\geq c'$. Then we have
\begin{equation*}
\fl N(x)=\lim_{T\to\infty}\left(\frac1{2\pi i}\int\limits_{c'-iT}^{c'+iT}
\zeta_{\Laplace}(s)x^s\frac{\dd s}s-
\sum_{2\pi|k|<T\log\lambda}\Res_{s=d_{\mathrm{S}}/2+2k\pi i/\log\lambda}
\frac{\zeta_{\Laplace}(s)x^s}s\right),
\end{equation*}
where we denote the real part of the poles of $R_w(\lambda^{-s})$ by $d_{\mathrm{S}}/2$
(the ``spectral dimension''). Now the limit of the integral
$\int_{c'-it}^{c'+iT}$ can be shown to exist for $T\to\infty$, which shows that
the limit
\begin{equation*}
\lim_{T\to\infty}\sum_{2\pi|k|<T\log\lambda}\Res_{s=d_{\mathrm{S}}/2+2k\pi i/\log\lambda}
\frac{\zeta_{\Laplace}(s)x^s}s
\end{equation*}
exists. This can be rewritten as $x^{d_{\mathrm{S}}/2}H(\log_\lambda x)$ for a periodic
continuous function $H$ given by its Fourier expansion. The limit of the
integral can be shown to be $\BigOh(x^{c'})$. Thus we have shown
\begin{equation}\label{eq:Nx-asymp}
N(x)=x^{d_{\mathrm{S}}/2}H(\log_\lambda x)+\BigOh(x^{c'}).
\end{equation}
Especially, the limit $\lim_{x\to\infty}x^{-d_{\mathrm{S}}/2}N(x)$ does not exist.

We remark here that there exists totally different approach to zeta functions
of fractals due to M.~Lapidus and his collaborators
\cite{He_Lapidus1997:generalized_minkowski_content,
Lapidus_Frankenhuysen1999:complex_dimensions_fractal,
Lapidus_Frankenhuysen2000:fractal_geometry_number,
Lapidus_Frankenhuijsen2006:fractal_geometry_complex}. In this geometric
approach the volume of tubular neighbourhoods of the fractal $G$
\begin{equation*}
V_G(\eps)=\vol\left(\left\{x\in \R^n\mid d(x,G)<\eps\right\}\right)
\end{equation*}
is analysed. The asymptotic behaviour of $V_G(\eps)$ for $\eps\to0$ gives rise
to the definition of a zeta function. In this geometric context the complex
solutions of \eref{eq:dimension} occur as poles of the zeta function; they are
called the ``complex dimensions'' of the fractal $G$ in this context. This
approach is motivated by the definition of Minkowski content, which itself
turned out to be too restrictive to measure (most of the) self-similar
fractals.

\subsection{Trace of the heat kernel}\label{sec:trace-heat-kernel}
The diffusion semi group $A_t$ introduced in
Section~\ref{sec:laplace-operator-inf} can be given in terms of the heat kernel
$p_t(x,y)$. As opposed to the situation in the euclidean case, the knowledge on
behaviour of the heat kernel for $t\to0$ is by far not as precise. Kumagai
\cite{Kumagai1993:estimates_brownian} proved the following lower
and upper estimates of the form
\begin{eqnarray}\label{eq:heat-estim-low}
t^{-\frac{d_{\mathrm{S}}}2}\exp
\left(-c_1\left(\frac{d(x,y)^{d_w}}t\right)^{\frac1{d_w-1}}\right)
\lesssim p_t(x,y)\\
 p_t(x,y)\lesssim
t^{-\frac{d_{\mathrm{S}}}2}\exp
\left(-c_2\left(\frac{d(x,y)^{d_w}}t\right)^{\frac1{d_w-1}}\right)
\label{eq:heat-estim-up}
\end{eqnarray}
where $d_{\mathrm{S}}$ and $d_w$ are the spectral and the walk
dimension of the fractal, respectively. These dimensions are related
to the Hausdorff dimension $d_f$ of the fractal via the so called
\emph{Einstein relation} $d_{\mathrm{S}}d_w=2d_f$. In the fractal case
usually $d_w>2$, as opposed to the euclidean case, where $d_w=2$,
which implies $d_{\mathrm{S}}=d_f$.

Only recently, a conjecture by Barlow and Perkins
\cite{Barlow_Perkins1988:brownian_motion_sierpinski} could be proved
by Kajino \cite{Kajino2011:on-diagonal_oscillations}. Namely that for
any $x\in K$ the limit
\begin{equation*}
\lim_{t\to0}t^{\frac{d_{\mathrm{S}}}2}p_t(x,x)
\end{equation*}
does \emph{not} exist for a large class of self-similar fractals. This gives
an indication, why obtaining more precise information on the heat kernel than
the estimates \eref{eq:heat-estim-low} and \eref{eq:heat-estim-up} would be
very difficult.

Even if the precise behaviour of the heat kernel seems to be far out of reach,
the trace of the heat kernel 
\begin{equation*}
K(t)=\int_K p_t(x,x)\dd\mathcal{H}(x)
\end{equation*}
can still be analysed in some detail under the
assumption of spectral decimation. The reason for this are the two relations
between $K(t)$, $N(x)$, and $\zeta_{\Laplace}$
\begin{eqnarray}\label{eq:P-N}
K(t)=\int_0^\infty e^{-tx}\dd N(x)=t\int_0^\infty N(x)e^{-xt}\dd x\\
\zeta_{\Laplace}(s)=\frac1{\Gamma(s)}\int_0^\infty K(t)t^{s-1}\dd t
\label{eq:P-zeta}.
\end{eqnarray}
The first expresses $K(t)$ as the Laplace transform of $N(x)$, the second gives
$\zeta_{\Laplace}$ as Mellin transform of $K(t)$.

Given the precise knowledge on the zeta function obtained in
Section~\ref{sec:spectr-zeta-funct}, we can use the Mellin inversion formula to
derive asymptotic information on $K(t)$ for $t\to0$. We have for
$c>\frac{d_{\mathrm{S}}}2$ 
\begin{equation}\label{eq:P-Mellin}
K(t)=\frac1{2\pi i}
\int\limits_{c-i\infty}^{c+i\infty}\zeta_{\Laplace}(s)\Gamma(s)t^{-s}\dd s,
\end{equation}
where integration is along the vertical line $\Re s=c$. We notice that the
Gamma function decays exponentially along vertical lines, whereas Dirichlet
series grow at most polynomially by the general theory of Dirichlet series
(cf.~\cite{Hardy_Riesz1964:general_theory_dirichlet}). Thus convergence of the
integral is guaranteed.

Shifting the line of integration in \eref{eq:P-Mellin} to $\Re s=-M-\frac12$
($M\in\mathbb{N}$) and taking the poles of the integrand into account, we
obtain for $t\to0$
\begin{equation}\label{eq:P-asymp}
K(t)=t^{-\frac{d_{\mathrm{S}}}2}H(\log_\lambda t)+\sum_j t^{-\alpha_j}H_j(\log_\lambda t)+
\BigOh(t^{M+\frac12}).
\end{equation}
Here, $H$ and $H_j$ are periodic continuous functions of period $1$, whose
Fourier coefficients are given as the residues of
$\zeta_{\Laplace}(s)\Gamma(s)$ in the poles on the lines $\Re s=\frac{d_{\mathrm{S}}}2$ or
$\Re s=\alpha_j$ respectively. By the strong decay of the Gamma function, these
functions are even real analytic. The values $\alpha_j+2k\pi i/\log\lambda$
come as the poles of the functions $R_w(\lambda^{-s})$. This shows that there
are only finitely many $\alpha_j$. Since $M$ can be made arbitrarily large, the
error term decays faster than any positive power of $t$ for $t\to0$.

The existence of complex poles of the zeta function thus implies the presence
of periodically oscillating terms in the asymptotic behaviour of the trace of
the heat kernel for $t\to0$. This implies that the limit
$\lim_{t\to0}t^{d_{\mathrm{S}}/2}p_t(x,x)$ does not exist on a set of positive measure for
$x$ in accordance with the above mentioned result by
N.~Kajino\cite{Kajino2011:on-diagonal_oscillations}. The consequences of the
presence of complex poles of the zeta function to properties of the heat
kernel, the density of states, and the partition function as well as the
physical implications of the resulting fluctuating behaviour in these
quantities have been discussed in
\cite{Akkermans_Dunne_Teplyaev2009:complex_dimensions,
  Akkermans_Dunne_Teplyaev2010:thermodynamics_fractals}.

\subsubsection{Casimir energy on fractals}\label{sec:casimir}
As an application of the spectral zeta function and its properties we compute
the Casimir energy of a fractal. We follow the lines of the exposition in
\cite{Kirsten2002:spectral_functions_mathematics} and refer to this book for
further details.

Consider the differential operator
$P=-\frac{\partial^2}{\partial\tau^2}+\Laplace$ on $(\R/\frac1\beta\Z)\times
G$. As usual the parameter $\beta=1/kT$. The zeta function of the operator $P$
is then given by
\begin{equation*}
\zeta_P(s)=\frac1{\Gamma(s)}\int_0^\infty K(t)
\sum_{n\in\Z}e^{-\frac{4\pi^2n^2}{\beta^2}t}t^{s-1}\,\dd t.
\end{equation*}
Using the theta function relation
\begin{equation*}
\sum_{n\in\Z}e^{-\frac{4\pi^2n^2}{\beta^2}t}=
\frac{\beta}{2\sqrt{\pi t}}\sum_{n\in\Z}e^{-\frac{\beta^2n^2}{4t}}
\end{equation*}
we obtain
\begin{equation}\label{eq:zetaP}
\fl \zeta_P(s)=\frac\beta{2\sqrt\pi\Gamma(s)}\Gamma\left(s-\frac12\right)
\zeta_{\Laplace}\left(s-\frac12\right)+
\frac\beta{\sqrt\pi\Gamma(s)}\int_0^\infty K(t)
\sum_{n=1}^\infty e^{-\frac{\beta^2n^2}{4t}} t^{s-\frac32}\,\dd t.
\end{equation}
The derivative $\zeta'_P(0)$ can be interpreted as $-\log\det P$, the logarithm
of the (regularised) determinant of $P$.

For the rescaled operator $P/\mu^2$ we have
\begin{equation*}
\zeta'_{P/\mu^2}(0)=\zeta'_P(0)+\zeta_P(0)\ln\mu^2.
\end{equation*}
From \eref{eq:zetaP} we compute
\begin{equation*}
\zeta_P(0)=-\beta\Res_{s=-\frac12}\zeta_{\Laplace}(s),
\end{equation*}
which vanishes in the case of the Sierpi\'nski gasket, because
$\zeta_{\Laplace}$ does not have a pole at $-\frac12$.
Furthermore, we get
\begin{equation*}
\zeta'_P(0)=-\beta\zeta_{\Laplace}\left(-\frac12\right)+
\frac\beta{\sqrt\pi}\sum_{n=1}^\infty\sum_{j=1}^\infty
\int_0^\infty e^{-\frac{\beta^2n^2}{4t}-\lambda_jt}t^{-\frac32}\,\dd t.
\end{equation*}
The integral and the summation over $n$ can be evaluated explicitly, which
finally gives
\begin{equation*}
\zeta'_P(0)=-\beta\zeta_{\Laplace}\left(-\frac12\right)-
2\sum_{j=1}^\infty\ln\left(1-e^{-\beta\sqrt{\lambda_j}}\right).
\end{equation*}
The case of the Sierpi\'nski gasket and other self-similar fractal shows an
important contrast to the manifold case, where the zeta function (generically)
has a pole at $-\frac12$ (see Section~\ref{sec:spectral_manifolds}). In the
manifold case the value $\zeta_{\Laplace}(-\frac12)$ has to be replaced by the
\emph{finite part} of $\zeta_{\Laplace}(s)$ at $s=-\frac12$, the function minus
the principal part at the polar singularity
(cf.~\cite{Kirsten2002:spectral_functions_mathematics}).

Now the energy of the system is given by
\begin{equation*}
E=-\frac12\frac{\partial}{\partial\beta}\zeta'_{P/\mu^2}(0)=
\frac12\zeta_{\Laplace}\left(-\frac12\right)+
\sum_{j=1}^\infty\frac{\sqrt{\lambda_j}}{e^{\beta\sqrt{\lambda_j}}-1}.
\end{equation*}
Letting $\beta\to\infty$, which is equivalent to letting temperature tend to
$0$, gives the Casimir energy
\begin{equation*}
E_{\mathrm{Cas}}=\frac12\zeta_{\Laplace}\left(-\frac12\right)
\end{equation*}
for fractals, whose zeta function has no pole at $-\frac12$.

\subsubsection{Numerical computations}\label{sec:numerical}
We will now describe the numerical computation of the values
$\zeta_{\Phi,w}(-1/2)$, which are needed for the computation of
$\zeta_{\Laplace}(-1/2)$ in the case of the Sierpi\'nski gasket. This
fractal admits spectral decimation with the polynomial $p(x)=x(x+5)$
(cf.~\cite{Rammal_Toulouse1983:random_walks_fractal,
  Fukushima_Shima1992:spectral_analysis_sierpinski}). The
corresponding Poincar\'e function is then given by the unique
holomorphic solution of the equation
\begin{equation}\label{eq:sierp-poinc}
\Phi(5z)=\Phi(z)(\Phi(z)+5),\quad \Phi(0)=0,\quad \Phi'(0)=1.
\end{equation}

We use the expression for $\zeta_{\Laplace}(s)$ for Dirichlet boundary
conditions derived in
\cite[Section~7]{Derfel_Grabner_Vogl2008:zeta_function_laplacian}
\begin{equation*}
\fl\zeta^{\mathrm{D}}_{\Laplace}(s)=5^{-s}\zeta_{\Phi,-2}(s)+
\frac3{5^s(5^s-1)(5^s-3)}\zeta_{\Phi,-3}(s)+
\frac{2\cdot5^s-5}{(5^s-1)(5^s-3)}\zeta_{\Phi,-5}(s).
\end{equation*}
Combining this with \eref{eq:cont-zeta} and inserting $s=-\frac12$ gives
\begin{equation}\label{eq:zeta1/2}\eqalign{
\fl\zeta^{\mathrm{D}}_{\Laplace}\left(-\frac12\right)=
\frac1{2\pi(\sqrt5-2)}\Biggl[
\sqrt5\int_0^\infty\log\Psi_{-2}(x)x^{-\frac12}\,\dd x\\
\!\!\!\!\!\!\!\!\!\!\!\!\!\!\!\!\!\!\!\!\!\!\!
+\frac{155+135\sqrt5}{88}\int_0^\infty\log\Psi_{-3}(x)x^{-\frac12}\,\dd x
-\frac{245+47\sqrt5}{124}\int_0^\infty\log\Psi_{-5}(x)x^{-\frac12}\,\dd x\Biggr].}
\end{equation}

From the general information on the asymptotic behaviour of Poincar\'e
functions derived in
\cite{Derfel_Grabner_Vogl2007:asymptotics_poincare_functions,
Derfel_Grabner_Vogl2008:complex_asymptotics_julia}
we obtain estimates of the form
\begin{equation*}
\exp(C_1x^\rho)\leq\Phi(x)\leq\exp(C_2x^\rho)
\end{equation*}
for positive constants $C_1$ and $C_2$, $\rho=\log_52$ and valid for $x\geq
x_0$. Such estimates can be proved easily by first showing the estimate for an
interval of the form $[x_0,5x_0]$ and then extending it by using the functional
equation \eref{eq:sierp-poinc}. For instance, we used $C_1=1$, $C_2=1.08$, and
$x_0=10$. Then the functions $\log\Psi_w(x)$ tend to $0$ like
$\exp(-C_2x^\rho)$ for $x\to\infty$. Given a precision goal $\eps>0$, we choose
$T$ so large that
\begin{equation*}
\int_T^\infty \exp(-C_2x^\rho)x^{-\frac12}\,\dd x<\eps.
\end{equation*}
Then the improper integrals in \eref{eq:zeta1/2} can be replaced by
$\int_0^T$. The functions $\Psi_w(x)$ can be computed to high precision by
using the power series representation of $\Phi(x)$ for $|x|\leq1$ and the
functional equation \eref{eq:sierp-poinc} to obtain
\begin{equation*}
\Phi(x)=p^{(k+1)}(\Phi(x/5^{k+1}))
\end{equation*}
for $5^k<|x|\leq 5^{k+1}$. This allows the computation of the remaining finite
integrals up to precision $\eps$. We obtained
\begin{equation*}
E^{\mathrm{D}}_{\mathrm{Cas}}=0.5474693544\ldots
\end{equation*}
for the Casimir energy of the two-dimensional Sierpi\'nski gasket with
Dirichlet boundary conditions.

Similarly, we have for Neumann boundary conditions
(cf.~\cite{Derfel_Grabner_Vogl2008:zeta_function_laplacian})
\begin{equation*}
\zeta^{\mathrm{N}}_{\Laplace}(s)=\frac1{(5^s-1)(5^s-3)}
\left((2\cdot5^s-5)\zeta_{\Phi,-3}(s)+5^s\zeta_{\Phi,-5}(s)\right).
\end{equation*}
This gives by the same numerical estimates as before
\begin{equation*}
E^{\mathrm{N}}_{\mathrm{Cas}}=2.134394089264\ldots
\end{equation*}
for the Casimir energy two-dimensional Sierpi\'nski gasket with
Neumann boundary conditions.

\subsection{Self-similarity and the renewal
  equation}\label{sec:self-simil-renew}
Let $K$ be a post-critically finite self similar fractal with self-similar
structure $(K,\Sigma,(F_i)_{i=1}^m)$ with contraction ratios $\alpha_i$ and
Hausdorff-dimension $\rho$. Let $K$ be equipped with a Dirichlet form with
parameters $r_i$ and $\lambda$ as described in
Section~\ref{sec:lapl-oper-dirichlet}. Denote by $\ND(x)$ and $\NN(x)$ the
eigenvalue counting functions for the Laplace operator under Dirichlet and
Neumann boundary conditions, respectively. Then the following crucial fact was
observed in \cite{Kigami_Lapidus1993:weyl's_problem_spectral}
\begin{equation}\label{eq:Nrenewal}
\sum_{i=1}^m \ND\left(\frac{r_i\alpha_i^\rho}\lambda x\right)\leq
\ND(x)\leq\NN(x)\leq
\sum_{i=1}^m \NN\left(\frac{r_i\alpha_i^\rho}\lambda x\right),
\end{equation}
where $\alpha_i$ are the contraction ratios introduced in
Section~\ref{sec:fractals}, $\rho$ is the Hausdorff dimension of $K$, and the
$r_i$ are the weights used for the construction of the Dirichlet form
$\mathcal{E}$ in Section~\ref{sec:lapl-oper-dirichlet}. Furthermore, the
inequalities
\begin{equation*}
\ND(x)\leq\NN(x)\leq\ND(x)+\#(V_0)
\end{equation*}
hold, where $V_0$ is defined in Definition~\ref{def:vm}. Setting
\begin{equation*}
\gamma_i=\left(\frac{r_i\alpha_i^\rho}\lambda\right)^{1/2},
\end{equation*}
we end up with the equation
\begin{equation}\label{eq:renewal}
\ND(x)=\sum_{i=1}^m\ND(\gamma_i^2x)+g(x),
\end{equation}
where $g(x)$ is defined as the difference between the left hand side and the
sum on the right hand side, and remains bounded by \eref{eq:Nrenewal}. A
similar equation holds for $\NN(x)$.

Equation \eref{eq:renewal} can be transformed into the classical renewal
equation occurring in probability theory
(cf.~\cite{Feller1966:introduction_probability_theory_ii}). The asymptotic
behaviour of the solutions of this equation is described by the following
theorem (stated as in \cite{Kigami2001:analysis_fractals}).
\begin{theorem}[``Renewal Theorem'']\label{thm:renewal}
Let $t^*>0$ and $f$ be a measurable function on $\R$ such that $f(t)=0$ for
$t<t^*$. If $f$ satisfies the renewal equation
\begin{equation*}
f(t)=\sum_{j=1}^Np_jf(t-\alpha_j)+u(t),
\end{equation*}
where $\alpha_1,\ldots,\alpha_N$ are positive real numbers, $p_j>0$ for
$j=1,\ldots,N$ and $\sum_{j=1}^Np_j=1$. Assume that $u$ is non-negative and
directly Riemann integrable on $\R$ with $u(t)=0$ for $t<t^*$. Then the
following conclusions hold:
\begin{enumerate}
\item arithmetic (or lattice) case: if the group generated by the
  $\alpha_j$ is discrete (i.e. there exist a $T>0$ and integers $m_j$
  with greatest common divisor $1$ such that $\alpha_j=Tm_j$; all the
  ratios $\alpha_i/\alpha_j$ are then rational), then
  $\lim_{t\to\infty}|f(t)-G(t)|=0$, where the $T$-periodic function
  $G$ is given by
  \begin{equation*}
    G(t)=\left(\sum_{j=1}^Np_jm_j\right)^{-1}\sum_{k=-\infty}^\infty
    u(t+kT).
  \end{equation*}
\item non-arithmetic (non-lattice) case: if the group generated by the
  $\alpha_j$ is dense in $\R$ (i.e. at least one of the ratios
  $\alpha_i/\alpha_j$ is irrational), then
  \begin{equation*}
    \lim_{t\to\infty}f(t)=\left(\sum_{j=1}^Np_j\alpha_j\right)^{-1}
\int_{-\infty}^\infty u(t)\,\dd t.
  \end{equation*}
\end{enumerate}
\end{theorem}
\begin{corollary}\label{cor:asymp}
Let $f$ be a solution of the equation
\begin{equation*}
f(x)=\sum_{i=1}^m f(\gamma_i^2x)+g(x),
\end{equation*}
where $0<\gamma_i<1$ and $g$ is a bounded function. Let $d_{\mathrm{S}}$ be the unique
positive solution of the equation
\begin{equation*}
\sum_{i=1}^m\gamma_i^{d_{\mathrm{S}}}=1,
\end{equation*}
then the following assertions hold:
\begin{enumerate}
\item arithmetic (or lattice) case: if the group generated by the values
  $\log\gamma_j$ is discrete, generated by $T>0$, then
  \begin{equation*}
    f(x)=x^{d_{\mathrm{S}}/2}\left(G((\log x)/2)+o(1)\right)
  \end{equation*}
for a periodic function $G$ of period $T$.
\item non-arithmetic (or non-lattice) case: if the group generated by the
  values $\log\gamma_j$ is dense, then
  \begin{equation*}
   \lim_{x\to\infty}f(x)x^{-d_{\mathrm{S}}/2}
  \end{equation*}
exists.
\end{enumerate}
\end{corollary}
The corollary is an immediate consequence of the Theorem by setting
$f(e^t)=e^{d_{\mathrm{S}}t/2}F(t)$ and applying the Theorem to $F$ and $g(e^t)e^{-d_{\mathrm{S}}t/2}$.

Summing up, we have the following theorem.
\begin{theorem}%
[{\cite[Theorem~2.4]{Kigami_Lapidus1993:weyl's_problem_spectral}}]
\label{thm:ev-asymp}
  Let $(K,\Sigma,(F_i)_{i=1}^m)$ be a self similar structure with
  contraction ratios $\alpha_i$ and Hausdorff-dimension $\rho$. Let
  $K$ be equipped with a Dirichlet form with parameters $r_i$ and
  $\lambda$ as described in Section~\ref{sec:lapl-oper-dirichlet}. Let
  $d_{\mathrm{S}}$ be the unique positive solution of the equation
  \begin{equation*}
    \sum_{i=1}^m\left(\frac{r_i\alpha_i^\rho}\lambda\right)^{d_{\mathrm{S}}/2}=1,
  \end{equation*}
  the ``spectral dimension'' of the harmonic structure of $K$.  Then the
  following assertions hold:
\begin{enumerate}
\item lattice case: if the group generated by the values
  $\log(r_i\alpha_i^\rho/\lambda)$ is dicrete, then
  \begin{equation*}
    \ND(x)=x^{d_{\mathrm{S}}/2}\left(G((\log(x))/2)+o(1)\right)
  \end{equation*}
for a periodic function $G$, which is non-constant in general.
\item non-lattice case: if the group generated by the values
  $\log(r_i\alpha_i^\rho/\lambda)$ is dense, then
  \begin{equation*}
    \lim_{x\to\infty}\ND(x)x^{-d_{\mathrm{S}}/2}
  \end{equation*}
exists.
\end{enumerate}
The behaviour of $\NN(x)$ for $x\to\infty$ is the same.
\end{theorem}

\begin{remark}\label{rem:levitin_vassiliev}
Similar ideas are used in \cite{Levitin_Vassiliev1996:spectral_asymptotics} to
study the asymptotic expansion of the eigenvalue counting function of the
Laplace operator on an open set $\mathfrak{G}$, which is formed from an open
set $G_0$ with smooth boundary as
\begin{equation*}
\mathfrak{G}=\bigcup_{n=0}^\infty G_n
\end{equation*}
with
\begin{equation*}
G_{n+1}=\bigcup_{i=1}^m F_i(G_n),
\end{equation*}
with similitudes $F_i$ and all the unions assumed to be disjoint. 
\end{remark}
\begin{remark}\label{rem:kajino}
  Very recently, N.~Kajino \cite{Kajino2010:spectral_laplacians} could
  extend this approach to non-p.~c.~f. fractals such as the
  Sierpi\'nski carpets. He observed that in order to have an
  inequality of the form \eref{eq:Nrenewal} it suffices to have the
  corresponding self-similarity property
\begin{equation*}
\mathcal{E}(f,f)=\lambda\sum_{i=1}^mr_i^{-1}\mathcal{E}(f\circ F_i,f\circ F_i)
\end{equation*}
of the underlying Dirichlet form. This self-similarity together with
symmetries characterises the Dirichlet form on the Sierpi\'nski carpet
uniquely, as was shown in
\cite{Barlow_Bass_Kumagai+2010:uniqueness_of_brownian}
\end{remark}

%% file: Self-similarity.tex
\section{Exploiting self-similarity}\label{sec:expl-self-simil}

\subsection{Decimation invariance and  equations with rescaling}
\label{sec:decim-invar-equat}
As has been described in Section~\ref{sec:random-walks-graphs},
diffusion processes on  decimation invariant finitely ramified
fractals, belonging to a certain class,
may be defined as  limits of  discrete simple symmetric
nearest-neighbour random walks on the associated approximating graphs.
The analysis is based on a renormalisation type argument, involving
self-similarity and decimation invariance .

Consider, for example, the Sierpi\'nski gasket $\Gamma$. It can be approximated by
``Sierpi\'nski lattice'' graphs $G_n$.  Let $X_{t}^{(n)}$ be a simple random
walk on $G_n$. Then, according to Goldstein
\cite{Goldstein1987:random_walks_diffusion} and Kusuoka
\cite{Kusuoka1987:diffusion_process_fractal},
 \begin{equation*}
 2^{-n}X_{[5^{n}t]}^{(n)} \Longrightarrow X_{t}
 \end{equation*}
 as $n \rightarrow \infty$, where $X_{t}$ is a diffusion process on $\Gamma$.
 ( Here, $\Longrightarrow$ means convergence in distribution and $[.]$ is the
 integer part function.)

 $X_{t}$ is a Markov process with continuous sample paths (in fact, a Feller
 process), which is itself invariant under the rescaling $x \rightarrow 2x,\;
 t \rightarrow 5t$.
  
 As has been explained in Section~\ref{sec:random-walks-graphs},
 equation~\eref{eq:green} in this instance has the form
 \begin{equation}\label{eq:rational}
  \Phi(\lambda z)= R(\Phi(z)), 
 \end{equation}
 where $\lambda=5;\: R(z)= \psi(z)=\frac{z^2}{4-3z}$ is a rational function of
 $z$.

Denoting here $\Phi(z)=1/\Psi(z)$ we have also
\begin{equation} \label{eq:polynomial}
\Psi(\lambda z)=P(\Psi(z))
\end{equation}
where $\lambda=5;\ P(z)=4z^2 -3z$ is a polynomial.

Both \eref{eq:rational} and \eref{eq:polynomial} are examples of the {\em
   Poincar\'e equation} which will be discussed below.

\subsection{Functional equations in the  theory of branching processes }\label{sec:branch}
Iterative functional equations, and the Poincar\'e equation. in particular,
occur also in the context of
 \emph{branching processes} (cf.~\cite{Harris1963:theory_branching_processes}).

Here a probability
generating function
\begin{equation*}
q(z)=\sum_{n=0}^\infty p_nz^n
\end{equation*}
encodes the offspring distribution, where  $p_n\geq0$ is the probability
that an individual has $n$ offsprings in the next generation (note that
$q(1)=1$). The growth rate $\lambda=q'(1)$ determines whether the population is
increasing ($\lambda>1$) or dying out ($\lambda\leq1$). In the first case the
branching process is called \emph{super-critical}. The probability generating
function $q^{(n)}(z)$ ($n$-th iterate of $q$) encodes the distribution of the
size $X_n$ of the $n$-th generation under the offspring distribution $q$. In
the case of a super-critical branching process it is known that the random
variables $\lambda^{-n}X_n$ tend to a limiting random variable $X_\infty$. The
moment generating function of this random variable
\begin{equation*}
f(z)=\mathbb{E}e^{-zX_\infty}
\end{equation*}
satisfies the functional equation
(cf.~\cite{Harris1963:theory_branching_processes})
\begin{equation} \label{eq:branching1}
f(\lambda z)=q(f(z)),
\end{equation}
which is  yet  another example of the \emph{Poincar\'e equation}.

If $q(z)= \frac{1}{P(1/z)}$, where $P$ is polynomial, then
\eref{eq:branching1} coincides with \eref{eq:polynomial}.

In general, such a reduction to \eref{eq:polynomial} is possible, when $q$ is
conjugate to a polynomial by a M\"{o}bius transformation.

\subsection{Historical remarks on the   Poincar\'e  equation }\label{sec:historical-remarks}

In his seminal papers
\cite{Poincare1886:une_classe_etendue,Poincare1890:une_classe_nouvelle}.
H.~Poincar\'e has studied the equation
\begin{equation}\label{Eq 1}
f(\lambda z)= R(f(z)),\quad z \in \mathbf{C},
\end{equation}
where $R(z)$ is a rational function and $\lambda\in\mathbf{C}$.  He proved
that, if $R(0)=0$, $R'(0)=\lambda$, and $|\lambda|>1$, then there exists a
meromorphic or entire solution of \eref{Eq 1}.  After Poincar\'e, \eref{Eq 1}
is called \emph{the Poincar\'e equation} and solutions of \eref{Eq 1} are
called \emph{the Poinca\'re functions }. The next important step was made by Valiron
\cite{Valiron1923:lectures_on_general,Valiron1954:fonctions_analytiques}, who
investigated the case, where $R(z)=P(z)$ is a polynomial, i.e.
\begin{equation}\label{eq:poincare}
f(\lambda z)=P(f(z)),\quad z \in \mathbf{C},
\end{equation}
and obtained conditions for the existence of an entire solution $f(z)$.
Furthermore, he derived the following asymptotic formula for
$M(r)=\max_{|z|\leq r}|f(z)|$:
\begin{equation}\label{Eq 3}
\ln M(r)\sim r^{\rho}Q\left(\frac{\ln r}{\ln |\lambda|}\right),
\quad r\rightarrow \infty.
\end{equation}
Here $Q(z)$ is a $1$-periodic function bounded between two positive constants,
$\rho=\frac{\ln m}{\ln |\lambda|}$ and $m=\deg P(z)$. 

Various aspects of the
Poincare functions have been studied  in the papers
\cite{Eremenko_Levin1989:periodic_points_polynomials,
Eremenko_Levin1992:estimation_characteristic_exponents,
 Eremenko_Sodin1990:iterations_rational_functions,
Ishizaki_Yanagihara2005:borel_and_julia,
Levin1991:pommerenke's_inequality,
 Romanenko_Sharkovsky2000:long_time_properties,
Derfel_Grabner_Vogl2007:asymptotics_poincare_functions,
Derfel_Grabner_Vogl2008:complex_asymptotics_julia}. 

\subsection{Applications: diffusion on fractals}\label{sec:applications}

In addition to \eref{Eq 3}, in applications (see Sections
\ref{sec:spectr-zeta-funct} and \ref{sec:trace-heat-kernel}
above, in particular ) it is important to know asymptotics of entire solutions
$f(z)$ in certain angular regions and even on specific rays $re^{i \vartheta}$
of the complex plane.

In the following section we present some recent results in this direction 
(\cite{Derfel_Grabner_Vogl2007:asymptotics_poincare_functions,
       Derfel_Grabner_Vogl2008:complex_asymptotics_julia}).

     Note, that in the aforementioned applications ( diffusion on fractals )
     scaling factor $\lambda$ is \emph{real} and, also, polynomial $P(z)$ has
     \emph{real} coefficients.

     Nevertheless, we prefer to start our exposition in section
     \ref{sec:asymptotics-spirals-rays} below from the general case
     and turn to the \emph{real} case afterwards.  We think that,
     these general facts on the Poincar\'e equation are interesting in
     their own, and hope that, maybe, they will also find applications
     in the future.

\subsection{Asymptotics along spirals and asymptotics along
rays}\label{sec:asymptotics-spirals-rays}

Here we intend to describe further results of Valiron's type .  We
start our presentation from some of our recent results
\cite{Derfel_Grabner_Vogl2007:asymptotics_poincare_functions,
  Derfel_Grabner_Vogl2008:complex_asymptotics_julia}.

It turns out that asymptotic behaviour of \emph{the Poinca\'re functions }
heavily depends on the arithmetic nature of $\arg\lambda$, where $\lambda$ is
the scaling factor in \eref{eq:poincare}. Denote $\arg\lambda=2\pi\beta$ .

The following statement is true, for \emph{irrational} $\beta$
(\cite{Derfel_Grabner_Vogl2007:asymptotics_poincare_functions}) :

\begin{theorem}\label{theorem1_irrational}
 If $\arg\lambda=2\pi\beta$ and $\beta$ is irrational, then $f(z)$
is unbounded along any ray $\vartheta$. Moreover, if we denote $\varphi(z)=
\ln |f(z)|$ (where main branch of logarithm is taken)
then there exists a sequence $r_n \rightarrow \infty$, such that the limit
\begin{equation}\label{Eq 5}
\lim_{n \rightarrow \infty}\frac{\varphi(r_ne^{i\vartheta})}{r_n^{\rho}}=L
\end{equation}
exists and $L>0$.
\end{theorem}

As far as we know this phenomenon has not been mentioned in the literature
before.

On the other hand, if  $\beta$ is rational (and, in particular, if
$\beta=0$, i.e. $\lambda$ is real) $f(z)$  may be bounded on some rays
and even in whole sectors. Nevertheless, for rational $\beta$,
the limit \eref{Eq 5} still exists under some additional assumptions.
Denote $\beta=t/s$ and suppose that $t$, $s$ are relatively prime.
Put $q=\lambda^s$ (note that $ 1<q \in \mathbf{R}$).

 The following result  is true for 
 \emph{rational}  $\beta$ and, in particular, for  \emph{real} $\lambda$ 
 (\cite{Derfel_Grabner_Vogl2007:asymptotics_poincare_functions}):

\begin{theorem}\label{theorem1_rational}
 Suppose that either $|\lambda|>m^2$ or $ s>2\rho$.
Then $f(z)$ is unbounded on any ray, and one can find a geometric progression
$r_n=q^n r_0$, $(r_0>0)$, for which the limit \eref{Eq 5} exists and $L>0$.
\end{theorem}
The above results are based on the following theorem on asymptotics of $f(z)$
along spirals (geometric progressions) of the form $z_{n}=\lambda^{n}z_{0},$
where $z_{0}=r_{0}e^{i\theta_{0}} \in \mathbb{C}$ will be specified below.

Consider the Poincar\'e equation
\begin{equation}\label{eq3.1}
f(\lambda z)= P(f(z)),
\end{equation}
where
\begin{equation}\label{eq3.2}
P(z) = z^{m} + p_{m-1}z^{m-1}+\ldots + p_1 z + p_{0}.
\end{equation}

Denote
\begin{equation}\label{eq3.3}
K = \max \{|p_{0}|,\ldots |p_{m-1}| \}.
\end{equation}
Suppose that $z_{0}=r_{0}e^{i\theta_{0}} \in \mathbb{C}$ is such a
point, that
\begin{equation}\label{eq3.4}
|f(z_{0})| > \max \{e, 2mK\}.
\end{equation}

Let us use the following notations
\begin{equation}\label{eq3.5}
z_{n}= \lambda^{n}z_{0};  \quad \phi(z)= \log|f(z)|;\quad  \rho= \frac{\log m}{\log|\lambda|}
\end{equation}
in part already introduced earlier. (It follows from \eref{eq3.4} that
$\phi(z_{0})> 1$ )

\begin{theorem}\label{theorem1}
Suppose that \eref{eq3.4} is satisfied. Then  the  limit
along the spiral $z_{n}= \lambda^{n}z_{0}$
\begin{equation}\label{eq3.6}
\lim_{n\rightarrow \infty}\frac{\phi(z_{n})}{|z_{n}|^\rho}= L(z_0)
\end{equation}
exists and
\begin{equation}\label{eq3.7}
\frac{\phi(z_{0})- 3Km/(2|f(z_{0})|}{|z_{0}|^{\rho}}< L(z_0) <
\frac{\phi(z_{0})+ 3Km/(2|f(z_{0})|}{|z_{0}|^{\rho}}.
\end{equation}
In particular, in view of \eref{eq3.4}
\begin{equation}\label{eq3.8}
\frac{\phi(z_{0})-3/4}{|z_{0}|^{\rho}}< L <
\frac{\phi(z_{0})+3/4}{|z_{0}|^{\rho}}.
\end{equation}
Furthermore, $L$ is a continuous function on the domain
$\{z\in\C\mid |f(z)|>\max(e,2mK)\}$.
\end{theorem}

\begin{remark}
  Here we deal mainly with the asymptotics of entire solutions of
  \eref{eq:poincare} However some statements are valid for arbitrary solutions
  of \eref{eq:poincare}.  In particular, for validity \eref{eq3.6} and
  \eref{eq3.7} \emph{no assumptions on the smoothness are needed.}
\end{remark}
\subsection{Asymptotics and dynamics in the real case. 
Properties of the Julia set }
\label{sec:real_julia}

Further refinements are possible when $\lambda>1$ is real and
$P(z)=p_mz^m +  \ldots +p_1z+p_0$ is a polynomial with real
coefficients.

Without loss of generality we can assume also that:
\begin{equation*}
f(0)=P(0)=0; \quad P'(0)=\lambda>1\mbox{ and } f'(0)=1
\end{equation*}

For reader convenience, we recall now some basic notions from the iteration
theory and complex dynamics \cite{Beardon1991:iteration_rational_functions,
  Falconer2003:fractal_geometry}.

To make things shorter, we give here definitions which slightly differ from
standard ones, but are equivalent to them in the \emph {polynomial} case.  That
will suffice our needs, in what follows.
  
We will especially need the
component of $\infty$ of  Fatou set  $\F(P)$ given by
\begin{equation}\label{eq:Fatou-infty}
\F_\infty(P)=\left\{z\in\C\mid \lim_{n\to\infty}P^{(n)}(z)=\infty\right\},
\end{equation}
The filled Julia set is given by
\begin{equation}\label{eq:filled-Julia}
\K(p)=\left\{z\in\C\mid (P^{(n)}(z))_{n\in\N}\mbox{ is bounded}\right\}=
\C\setminus\F_\infty(P).
\end{equation}
Now, according to ~\cite{Falconer2003:fractal_geometry} one can define the
Julia set $\J(p)$, as
\begin{equation}\label{eq:boundary}
\partial\K(P)=\partial\F_\infty(P)=\J(P).
\end{equation}
In the case of polynomials this can be used as an equivalent definition of
the Julia set.

\medskip Let $f(z)$ be an entire solution of \eref{eq:poincare}. In contrast
with the previous section (where asymptotics of $|f(z)|$ is studied) here we
collect some results on the asymptotics of the solution $f(z)$ \emph{itself}
(in some angular regions of the complex plane).  Our presentation is based on
(\cite{Derfel_Grabner_Vogl2008:zeta_function_laplacian,
  Derfel_Grabner_Vogl2007:asymptotics_poincare_functions,
  Derfel_Grabner_Vogl2008:complex_asymptotics_julia}).

\begin{theorem}
[{\cite[Theorem~1]{Derfel_Grabner_Vogl2008:zeta_function_laplacian}}]
\label{thm1}
Let $f$ be an entire solution of the functional equation \eref{eq:poincare}.
Furthermore, suppose that $\mathcal{F}_\infty$, the Fatou component of $\infty$
of $P$, contains an angular region of the form
\begin{equation*}
W_\beta=\left\{z\in\mathbb{C}\setminus\{0\}\mid |\arg z|<\beta\right\}
\end{equation*}                                                                
for some $\beta>0$. Then for any $\eps>0$ and any                              
$M>0$ the asymptotic relation                                                  
\begin{equation}\label{eq:Phiasymp}                                            
f(z)=\exp\left(z^\rho Q\left(\frac{\log z}{\log\lambda}\right)+                
o\left(|z|^{-M}\right)\right)                                                  
\end{equation}                                                                 
holds uniformly for $z\in W_{\beta-\eps}$, where $Q$ is a periodic holomorphic 
function of period $1$ on the strip $\{w\in\mathbb{C}\mid|\Im                  
w|<\frac\beta{\log\lambda}\}$. The real part of $z^\rho Q(\frac{\log           
  z}{\log\lambda})$ is bounded between two positive constants; $Q$ takes real  
values on the real axis.                                                       
\end{theorem}
\begin{remark}
  Notice that the condition on the Fatou component $\mathcal{F}_\infty$ is used
in the proof of this theorem to ensure that $f(z)$ tends to infinity in the
angular region $W_\beta$. Therefore, this condition could be replaced by
\begin{equation*}
\lim_{z\to\infty}f(z)=\infty\mbox{ for }|\arg z|<\beta.
\end{equation*}
\end{remark}

Yet a stronger result can be derived under the additional assumption
that the Julia set $\J(P)$ of polynomial $P(z)$ is \emph {real} (see
Corollary~\ref{corr:real_julia}, below).
\begin{remark}
  The latter assumption on the reality of $\J(P)$ looks artificial at the first
  glance, but, in fact is typical in applications, related to the diffusion on
  fractals.  The rough explanation of the latter fact is the following:

  The zeros of the solution $f(z)$ for \eref{eq:poincare} are eigenvalues of
  the infinitesimal generator of the diffusion, i.e. the ``Laplacian on the
  fractal''. The latter operator is self-adjoint , and its eigenvalues are
  real. Therefore zeros of $f(z)$ are real and, finally, this implies that
  Julia set of polynomial $P(z)$ should be real. This motivates our special
  interest for Poincar\'e equations with polynomial $P(z)$, having real Julia
  set.
\end{remark}

\begin{corollary}%
[{\cite[Corollary~4.1]{Derfel_Grabner_Vogl2008:complex_asymptotics_julia}}]
\label{corr:real_julia}
  Assume that $P$ is a real polynomial such that $\J(P)$ is real and all
  coefficients $p_i$ ($i\geq2$) of $P$ are non-negative.  Then
  $\J(P)\subset\R^-\cup\{0\}$ and therefore
\begin{equation}\label{eq:simple-asymp}
f(z)\sim\exp\left(z^\rho Q\left(\frac{\log z}{\log\lambda}\right)\right)
\end{equation}
for $z\to\infty$ and $|\arg z|<\pi$. Here $Q$ is a periodic function of period
$1$ holomorphic in the strip given by 
$|\Im w|<\frac\pi{\log\lambda}$. Furthermore, for every $\eps>0$
 the real part of $z^\rho Q(\frac{\log z}{\log\lambda})$ is bounded between two
positive constants  for $|\arg z|\leq\pi-\eps$.
\end{corollary}

If $P(z)$ is a quadratic polynomial ( a case, arising in numerous applications)
it is possible to give an exact criterion for reality of   $\J(P)$:
\medskip
\begin{lemma}%
[{\cite[Lemma~6.7]{Derfel_Grabner_Vogl2007:asymptotics_poincare_functions}}]
\label{lemma:quadratic_julia}
  Let
\begin{equation}\label{Eq 7}
P(z)=az(z- \omega), \quad 0\neq \omega\in\mathbf{R}
\end{equation}
Then Julia set $\J(P)$  is real, if and only if the following condition
is fulfilled
\begin{equation}\label{Eq 8}
a|\omega|\ge
\left\{
\begin{array}{lr}
2,&\omega>0\\
4,&\omega<0
\end{array}
\right.
\end{equation}
\end{lemma}
\medskip


It would be interesting to find some constructive  \emph{conditions for reality}
of   $\J(P)$  in terms of coefficients of $P$
(at least for polynomial of small order in the beginning).

Note also, that in the above Lemma $|P'(0)|= a|\omega |$ and $m=\deg P(z)=2 $.

Therefore \eref{Eq 8} can be rewritten in the form:
\begin{equation}\label{Eq 9}
 |P'(0)|\ge
\left\{
\begin{array}{lr}
m,&\omega>0\\
m^2,&\omega<0
\end{array}
\right.
\end{equation}
It turns out that the \emph{necessity}  part of latter criterion  (in the form \eref{Eq 9} ) is valid for general polynomial $P$, 
with \emph{real} Julia set $\J(P)$. 
Namely, the following result  of Pommerenke--Levin--Eremenko--Yoccoz type
(on inequalities for multipliers ) is true:

\begin{theorem}
[{\cite[Theorem ~4.1]{Derfel_Grabner_Vogl2008:complex_asymptotics_julia}}]
\label{thm:pommerenke_levin_eremenko}
Let $P$ be a polynomial of degree $m>1$ with real Julia set $\J(P)$. Then for
any fixed point $\xi$ of $P$ with $\min \J(P)<\xi<\max \J(P)$ we have
$|P'(\xi)|\geq m$. Furthermore, $|P'(\min\J(P))|\geq m^2$ and $|P'(\max
\J(P))|\geq m^2$. Equality in one of these inequalities implies that $P$ is
linearly conjugate to the Chebyshev polynomial $T_m$ of degree $m$.
\end{theorem}
\begin{remark}
This theorem can be compared to 
\cite{Levin1991:pommerenke's_inequality,
Pommerenke1986:conformal_mapping_iteration,
Eremenko_Levin1992:estimation_characteristic_exponents,
Buff2003:bieberbach_conjecture_dynamics}
 where inequalities (of the opposite direction ) for the
 multipliers of $P$ with \emph{connected} Julia sets were derived. 
\end{remark}
In the foregoing we dealt mainly with the asymptotics of solutions for the
polynomial Poincar\'e equation, only.

The study of the asymptotic behaviour of meromorphic solutions of the general
Poincar\'e equation

\begin{equation*}
f(\lambda z)= R(f(z)),\quad z \in \mathbf{C},
\end{equation*}
where $R(z)$ is rational function (rather than polynomial $P(z)$) still, to
large extent, remains an open challenge.


%% file: Physics_Review.bbl
\begin{thebibliography}{100}

\bibitem{Mandelbrot1982:fractal_geometry_nature}
B.~B. Mandelbrot.
\newblock {\em The fractal geometry if nature}.
\newblock Freeman, San Fransisco, 1982.

\bibitem{Alexander_Orbach1982:density_states_fractals}
S.~Alexander and R.~Orbach.
\newblock Density of states on fractals: ``fractons''.
\newblock {\em J. Physique Lettres}, 43:L625--L631, 1982.

\bibitem{Rammal_Toulouse1983:random_walks_fractal}
R.~Rammal and G.~Toulouse.
\newblock Random walks on fractal structures and percolation clusters.
\newblock {\em J. Physique Lettres}, 44:L13--L22, 1983.

\bibitem{Havlin_Ben-Avraham1987:diffusion_disordered_media}
S.~Havlin and D.~Ben-Avraham.
\newblock Diffusion on disordered media.
\newblock {\em Adv. Phys.}, 36:696--798, 1987.

\bibitem{Havlin_Ben-Avraham2000:diffusion_reactions_fractals}
D.~Ben-Avraham and S.~Havlin.
\newblock {\em Diffusion and reactions in fractals and disordered systems}.
\newblock Cambridge University Press, 2000.

\bibitem{Goldstein1987:random_walks_diffusion}
S.~Goldstein.
\newblock Random walks and diffusions on fractals.
\newblock In {\em Percolation theory and ergodic theory of infinite particle
  systems ({M}inneapolis, {M}inn., 1984--1985)}, volume~8 of {\em IMA Vol.
  Math. Appl.}, pages 121--129. Springer, New York, 1987.

\bibitem{Kusuoka1987:diffusion_process_fractal}
S.~Kusuoka.
\newblock A diffusion process on a fractal.
\newblock In {\em Probabilistic methods in mathematical physics
  ({K}atata/{K}yoto, 1985)}, pages 251--274. Academic Press, Boston, MA, 1987.

\bibitem{Barlow_Perkins1988:brownian_motion_sierpinski}
M.~T. Barlow and E.~A. Perkins.
\newblock Brownian motion on the {S}ierpi\'nski gasket.
\newblock {\em Probab. Theory Relat. Fields}, 79:543--623, 1988.

\bibitem{Harris1963:theory_branching_processes}
T.~E. Harris.
\newblock {\em The {T}heory of {B}ranching {P}rocesses}.
\newblock Springer, Berlin, New York, 1963.

\bibitem{Lindstroem1990:brownian_motion_nested}
T.~Lindstr{\o}m.
\newblock {\em Brownian {M}otion on {N}ested {F}ractals}, volume 420 of {\em
  Mem. Amer. Math. Soc.}
\newblock Amer. Math. Soc., 1990.

\bibitem{Barlow1998:diffusion_fractals}
M.~T. Barlow.
\newblock Diffusion on fractals.
\newblock In P.~Bernard, editor, {\em Lectures on Probability Theory and
  Statistics}, volume 1690 of {\em Lecture Notes in Mathematics}, pages 1--121.
  Springer, Berlin, 1998.

\bibitem{deGennes1976:percolation_concept}
P.~G. de~Gennes.
\newblock La percolation: un concept unifacateur.
\newblock {\em La Recherche}, 7:912--927, 1976.

\bibitem{Kigami1989:harmonic_calculus_sierpinski}
J.~Kigami.
\newblock A harmonic calculus on the {S}ierpi\'nski spaces.
\newblock {\em Japan J. Appl. Math.}, 6:259--290, 1989.

\bibitem{Kigami1993:harmonic_calculus_pcf}
J.~Kigami.
\newblock Harmonic calculus on p.c.f.\ self-similar sets.
\newblock {\em Trans. Amer. Math. Soc.}, 335:721--755, 1993.

\bibitem{Kigami1998:distributions_localized_eigenvalues}
J.~Kigami.
\newblock Distributions of localized eigenvalues of {L}aplacians on post
  critically finite self-similar sets.
\newblock {\em J. Funct. Anal.}, 156:170--198, 1998.

\bibitem{Kigami2001:analysis_fractals}
J.~Kigami.
\newblock {\em Analysis on fractals}, volume 143 of {\em Cambridge Tracts in
  Mathematics}.
\newblock Cambridge University Press, Cambridge, 2001.

\bibitem{Fukushima_Shima1992:spectral_analysis_sierpinski}
M.~Fukushima and T.~Shima.
\newblock On the spectral analysis for the {S}ierpi\'nski gasket.
\newblock {\em J. of Potential Analysis}, 1:1--35, 1992.

\bibitem{Kusuoka1989:dirichlet_forms_random_matrices}
S.~Kusuoka.
\newblock Dirichlet forms on fractals and products of random matrices.
\newblock {\em Publ. Res. Inst. Math. Sci.}, 25:659--680, 1989.

\bibitem{Barlow_Bass1989:construction_brownian_motion}
M.~T. Barlow and R.~F. Bass.
\newblock The construction for {B}rownian motion on the {S}ierpi\'nski carpet.
\newblock {\em Ann. Inst. H. Poincar\'e}, 25:225--257, 1989.

\bibitem{Shima1993:eigenvalue_problem_laplacian}
T.~Shima.
\newblock The eigenvalue problem for the {L}aplacian on the {S}ierpi\'nski
  gasket.
\newblock In K.D. Elworthy and N.~Ikeda, editors, {\em Asymptotic problems in
  probability theory: stochastic models and diffusions on fractals
  (Sanda/Kyoto, 1990)}, volume 283 of {\em Pitman Res. Notes Math. Ser.}, pages
  279--288. Longman Sci. Tech., Harlow, 1993.

\bibitem{Shima1996:eigenvalue_problems_laplacians}
T.~Shima.
\newblock On eigenvalue problems for {L}aplacians on p.c.f. self-similar sets.
\newblock {\em Japan J. Indust. Appl. Math.}, 13:1--23, 1996.

\bibitem{Malozemov_Teplyaev2003:self_similarity_operators}
L.~Malozemov and A.~Teplyaev.
\newblock Self-similarity, operators and dynamics.
\newblock {\em Math. Phys. Anal. Geom.}, 6:201--218, 2003.

\bibitem{Strichartz2003:fractafolds_spectra}
R.~S. Strichartz.
\newblock Fractafolds based on the {S}ierpi\'nski gasket and their spectra.
\newblock {\em Trans. Amer. Math. Soc.}, 355:4019--4043, 2003.

\bibitem{Teplyaev2004:spectral_zeta_function}
A.~Teplyaev.
\newblock Spectral zeta function of symmetric fractals.
\newblock In C.~Bandt, U.~Mosco, and M.~Z{\"a}hle, editors, {\em Fractal
  geometry and stochastics III}, volume~57 of {\em Progr. Probab.}, pages
  245--262. Birkh\"auser, Basel, 2004.

\bibitem{Beardon1991:iteration_rational_functions}
A.~F. Beardon.
\newblock {\em Iteration of rational functions}.
\newblock Number 132 in Graduate Texts in Mathematics. Springer Verlag, 1991.

\bibitem{Milnor2006:dynamics_complex}
J.~Milnor.
\newblock {\em Dynamics in one complex variable}, volume 160 of {\em Annals of
  Mathematics Studies}.
\newblock Princeton University Press, Princeton, N. J, third edition, 2006.

\bibitem{Lapidus1992:spectral_fractal_geometry}
M.~L. Lapidus.
\newblock Spectral and fractal geometry: from the {W}eyl-{B}erry conjecture for
  the vibrations of fractal drums to the {R}iemann zeta-function.
\newblock In {\em Differential equations and mathematical physics
  ({B}irmingham, {AL}, 1990)}, volume 186 of {\em Math. Sci. Engrg.}, pages
  151--181. Academic Press, Boston, MA, 1992.

\bibitem{Lapidus_Frankenhuysen2000:fractal_geometry_number}
M.~L. Lapidus and M.~van Frankenhuysen.
\newblock {\em Fractal geometry and number theory, \emph{Complex dimensions of
  fractal strings and zeros of zeta functions}}.
\newblock Birkh{\"a}user Boston Inc., Boston, MA, 2000.

\bibitem{Lapidus_Frankenhuijsen2006:fractal_geometry_complex}
M.~L. Lapidus and M.~van Frankenhuysen.
\newblock {\em Fractal geometry, complex dimensions and zeta functions,
  \emph{Geometry and spectra of fractal strings}}.
\newblock Springer Monographs in Mathematics. Springer, New York, 2006.

\bibitem{Akkermans_Dunne_Teplyaev2009:complex_dimensions}
E.~Akkermans, G.~Dunne, and A.~Teplyaev.
\newblock Physical consequences of complex dimensions of fractals.
\newblock {\em Europhys. Lett. EPL}, 88:40007, 2009.
\newblock arXiv:0903.3681.

\bibitem{Akkermans_Dunne_Teplyaev2010:thermodynamics_fractals}
E.~Akkermans, G.~Dunne, and A.~Teplyaev.
\newblock Thermodynamics of photons on fractals.
\newblock {\em Phys. Rev. Lett.}, 105:230407, 2010.
\newblock arXiv:1010.1148.

\bibitem{Strichartz2006:differential_equations_fractals}
R.~S. Strichartz.
\newblock {\em Differential equations on fractals}.
\newblock Princeton University Press, Princeton, NJ, 2006.
\newblock A tutorial.

\bibitem{Triebel1997:fractals_and_spectra}
H.~Triebel.
\newblock {\em Fractals and spectra}, volume~91 of {\em Monographs in
  Mathematics}.
\newblock Birkh\"auser Verlag, Basel, 1997.
\newblock Related to Fourier analysis and function spaces.

\bibitem{Kirillov2009:fractals}
A.~A. Kirillov.
\newblock {\em A tale of two fractals (Russian)}.
\newblock MCCME, Moscow, 2009.
\newblock preliminary English version available at
  \url{http://www.math.upenn.edu/\~{}kirillov/MATH480-F07/tf.pdf}.

\bibitem{Kumagai2008:recent_developments_fractals}
T.~Kumagai.
\newblock Recent developments of analysis on fractals.
\newblock {\em Translations of American Math. Soc. series 2,}, 223:81--95,
  2008.

\bibitem{Teplyaev2010:diffusions_and_spectral}
A.~Teplyaev.
\newblock Diffusions and spectral analysis on fractals: an overview.
\newblock available at
  \url{http://www.math.uconn.edu/\~{}teplyaev/POMI-2010-teplyaev.pdf},
  September 2010.
\newblock talk given at ``Mathematics - XXI century'', conference for the 70th
  anniversary of the St.~Petersburg Department of V.~A.~Steklov Institute of
  Mathematics of the Russian Academy of Sciences.

\bibitem{Barlow_Bass1999:brownian_sierpinski_carpet}
M.~T. Barlow and R.~F. Bass.
\newblock Brownian motion and harmonic analysis on {S}ierpinski carpets.
\newblock {\em Canad. J. Math.}, 51(4):673--744, 1999.

\bibitem{Barlow_Bass_Kumagai+2010:uniqueness_of_brownian}
M.~T. Barlow, R.~F. Bass, T.~Kumagai, and A.~Teplyaev.
\newblock Uniqueness of {B}rownian motion on {S}ierpi\'nski carpets.
\newblock {\em J. Eur. Math. Soc. (JEMS)}, 12(3):655--701, 2010.

\bibitem{Barlow_Coulhon_Grigoryan2001:manifolds_graphs_with}
M.~Barlow, T.~Coulhon, and A.~Grigor'yan.
\newblock Manifolds and graphs with slow heat kernel decay.
\newblock {\em Invent. Math.}, 144:609--649, 2001.

\bibitem{Grigoryan2010:heat_kernels_metric}
A.~Grigor'yan.
\newblock Heat kernels on metric measure spaces with regular volume growth.
\newblock In {\em Handbook of geometric analysis, {N}o. 2}, volume~13 of {\em
  Adv. Lect. Math. (ALM)}, pages 1--60. Int. Press, Somerville, MA, 2010.

\bibitem{Grigoryan_Hu_Lau2010:comparison_inequalities}
A.~Grigor'yan, J.~Hu, and K.-S. Lau.
\newblock Comparison inequalities for heat semigroups and heat kernels on
  metric measure spaces.
\newblock {\em J. Funct. Anal.}, 259(10):2613--2641, 2010.

\bibitem{Barlow_Grigoryan2009:heat_kernel_upper}
M.~T. Barlow and T.~Grigor'yan, A.and~Kumagai.
\newblock Heat kernel upper bounds for jump processes and the first exit time.
\newblock {\em J. Reine Angew. Math.}, 626:135--157, 2009.

\bibitem{Naimark_Solomyak1994:eigenvalue_behaviour}
K.~Naimark and M.~Solomyak.
\newblock On the eigenvalue behaviour for a class of operators related to
  self-similar measures on {${\bf R}^d$}.
\newblock {\em C. R. Acad. Sci. Paris S{\'e}r. I Math.}, 319(8):837--842, 1994.

\bibitem{Naimark_Solomyak1995:eigenvalue_behaviour}
K.~Naimark and M.~Solomyak.
\newblock The eigenvalue behaviour for the boundary value problems related to
  self-similar measures on {${\bf R}^d$}.
\newblock {\em Math. Res. Lett.}, 2(3):279--298, 1995.

\bibitem{Naimark_Solomyak2001:eigenvalue_distribution_fractal}
K.~Naimark and M.~Solomyak.
\newblock Eigenvalue distribution of some fractal semi-elliptic differential
  operators: combinatorial approach.
\newblock {\em Integral Equations Operator Theory}, 40(4):495--506, 2001.

\bibitem{Solomyak_Verbitsky1995:spectral_problem_related}
M.~Solomyak and E.~Verbitsky.
\newblock On a spectral problem related to self-similar measures.
\newblock {\em Bull. London Math. Soc.}, 27(3):242--248, 1995.

\bibitem{Triebel2008:fractal_analysis_approach}
H.~Triebel.
\newblock Fractal analysis, an approach via function spaces.
\newblock In {\em Topics in mathematical analysis}, volume~3 of {\em Ser. Anal.
  Appl. Comput.}, pages 413--447. World Sci. Publ., Hackensack, NJ, 2008.

\bibitem{Hu_Lau_Ngai2006:laplace_operators_related}
J.~Hu, K.-S. Lau, and S.-M. Ngai.
\newblock Laplace operators related to self-similar measures on
  {$\mathbb{R}^d$}.
\newblock {\em J. Funct. Anal.}, 239(2):542--565, 2006.

\bibitem{Falconer1986:geometry_fractal_sets}
K.~J. Falconer.
\newblock {\em The {G}eometry of {F}ractal {S}ets}, volume~85 of {\em Cambridge
  Tracts in Mathematics}.
\newblock Cambridge University Press, Cambridge, 1986.

\bibitem{Hutchinson1981:fractals_self_similarity}
J.~E. Hutchinson.
\newblock Fractals and self-similarity.
\newblock {\em Indiana Univ. Math. J.}, 30:713--747, 1981.

\bibitem{Falconer1988:hausdorff_dimension_self}
K.~J. Falconer.
\newblock The {H}ausdorff dimension of self-affine fractals.
\newblock {\em Math. Proc. Cambridge Philos. Soc.}, 103:339--350, 1988.

\bibitem{Falconer1992:dimension_self_affine_ii}
K.~J. Falconer.
\newblock The dimension of self-affine fractals. {II}.
\newblock {\em Math. Proc. Cambridge Philos. Soc.}, 111:169--179, 1992.

\bibitem{Falconer1997:techniques_in_fractal}
K.~J. Falconer.
\newblock {\em Techniques in {F}ractal {G}eometry}.
\newblock John Wiley \& Sons Ltd., Chichester, 1997.

\bibitem{Hata1985:structure_self-similar}
M.~Hata.
\newblock On the structure of self-similar sets.
\newblock {\em Japan J. Appl. Math.}, 2(2):381--414, 1985.

\bibitem{Rosenberg1997:laplacian_riemannian_manifold}
S.~Rosenberg.
\newblock {\em The {L}aplacian on a {R}iemannian manifold}, volume~31 of {\em
  London Mathematical Society Student Texts}.
\newblock Cambridge University Press, Cambridge, 1997.

\bibitem{Federer1969:geometric_measure_theory}
H.~Federer.
\newblock {\em Geometric measure theory}.
\newblock Die Grundlehren der mathematischen Wissenschaften, Band 153.
  Springer-Verlag New York Inc., New York, 1969.

\bibitem{Fukushima_Oshima_Takeda2011:dirichlet_forms_markov}
M.~Fukushima, Y.~Oshima, and M.~Takeda.
\newblock {\em Dirichlet forms and symmetric {M}arkov processes}, volume~19 of
  {\em de Gruyter Studies in Mathematics}.
\newblock Walter de Gruyter \& Co., Berlin, extended edition, 2011.

\bibitem{Yosida1971:functional_analysis}
K.~Yosida.
\newblock {\em Functional {A}nalysis}.
\newblock Springer, Berlin, 1971.

\bibitem{Kroen2002:green_functions_self}
B.~Kr{\"o}n.
\newblock Green functions on self-similar graphs and bounds for the spectrum of
  the {L}aplacian.
\newblock {\em Ann. Inst. Fourier (Grenoble)}, 52:1875--1900, 2002.

\bibitem{Kroen_Teufl2004:asymptotics_transition_probabilities}
B.~Kr{\"o}n and E.~Teufl.
\newblock Asymptotics of the transition probabilities of the simple random walk
  on self-similar graphs.
\newblock {\em Trans. Amer. Math. Soc.}, 356:393--414, 2004.

\bibitem{Metz1995:hilberts_projective_metric}
V.~Metz.
\newblock Hilbert's projective metric on cones of {D}irichlet forms.
\newblock {\em J. Funct. Anal.}, 127(2):438--455, 1995.

\bibitem{Metz1996:renormalization_contracts}
V.~Metz.
\newblock Renormalization contracts on nested fractals.
\newblock {\em J. Reine Angew. Math.}, 480:161--175, 1996.

\bibitem{Metz2003:cone_diffusions}
V.~Metz.
\newblock The cone of diffusions on finitely ramified fractals.
\newblock {\em Nonlinear Anal.}, 55(6):723--738, 2003.

\bibitem{Sabot1997:existence_uniqueness_diffusion}
C.~Sabot.
\newblock Existence and uniqueness of diffusions on finitely ramified
  self-similar fractals.
\newblock {\em Ann. Sci. \'Ecole Norm. Sup. (4)}, 30(5):605--673, 1997.

\bibitem{Teufl2007:asymptotic_behaviour_analytic}
E.~Teufl.
\newblock On the asymptotic behaviour of analytic solutions of linear iterative
  functional equations.
\newblock {\em Aequationes Math.}, 73:18--55, 2007.

\bibitem{Grabner_Woess1997:functional_iterations_periodic}
P.~J. Grabner and W.~Woess.
\newblock Functional iterations and periodic oscillations for the simple random
  walk on the {S}ierpi\'nski graph.
\newblock {\em Stochastic Processes Appl.}, 69:127--138, 1997.

\bibitem{Grabner1997:functional_iterations_stopping}
P.~J. Grabner.
\newblock Functional iterations and stopping times for {B}rownian motion on the
  {S}ierpi\'nski gasket.
\newblock {\em Mathematika}, 44:374--400, 1997.

\bibitem{Schweitzer2006:diffusion_simply_doubly}
G.~Schweitzer.
\newblock {\em Diffusion on {S}imply and {D}oubly {S}ymmetric {F}ractal
  {G}raphs}.
\newblock PhD thesis, Graz University of Technology, 2006.

\bibitem{Derfel_Grabner_Vogl2008:zeta_function_laplacian}
G.~Derfel, P.~J. Grabner, and F.~Vogl.
\newblock The {Z}eta function of the {L}aplacian on certain fractals.
\newblock {\em Trans. Amer. Math. Soc.}, 360:881--897, 2008.

\bibitem{Barlow_Bass1999:graphical_sierpinski_carpet}
M.~T. Barlow and R.~F. Bass.
\newblock Random walks on graphical {S}ierpinski carpets.
\newblock In {\em Random walks and discrete potential theory ({C}ortona,
  1997)}, Sympos. Math., XXXIX, pages 26--55. Cambridge Univ. Press, Cambridge,
  1999.

\bibitem{Ben-Bassat_Strichartz_Teplyaev1999:not_domain_laplacian}
O.~Ben-Bassat, R.~S. Strichartz, and A.~Teplyaev.
\newblock What is not in the domain of the {L}aplacian on {S}ierpinski gasket
  type fractals.
\newblock {\em J. Funct. Anal.}, 166(2):197--217, 1999.

\bibitem{Fukushima1992:dirichlet_forms}
M.~Fukushima.
\newblock Dirichlet forms, diffusion processes and spectral dimensions for
  nested fractals.
\newblock In S.~Albeverio, J.~E. Fenstad, H.~Holden, and T.~Lindstr{\o}m,
  editors, {\em Ideas and methods in mathematical analysis, stochastics, and
  applications ({O}slo, 1988)}, pages 151--161. Cambridge Univ. Press,
  Cambridge, 1992.

\bibitem{Minakshisundaram_Pleijel1949:some_properties_eigenfunctions}
S.~Minakshisundaram and {\AA}.~Pleijel.
\newblock Some properties of the eigenfunctions of the {L}aplace-operator on
  {R}iemannian manifolds.
\newblock {\em Canadian J. Math.}, 1:242--256, 1949.

\bibitem{Berline_Getzler_Vergne1992:heat_kernels}
N.~Berline, E.~Getzler, and M.~Vergne.
\newblock {\em Heat kernels and {D}irac operators}, volume 298 of {\em
  Grundlehren der Mathematischen Wissenschaften}.
\newblock Springer-Verlag, Berlin, 1992.

\bibitem{Polterovich2000:heat_invariants}
I.~Polterovich.
\newblock Heat invariants of {R}iemannian manifolds.
\newblock {\em Israel J. Math.}, 119:239--252, 2000.

\bibitem{Kirsten2002:spectral_functions_mathematics}
K.~Kirsten.
\newblock {\em Spectral Functions in Mathematics and Physics}.
\newblock Chapman \& Hall/CRC, 2002.

\bibitem{Kac1966:hear_drum}
M.~Kac.
\newblock Can one hear the shape of a drum?
\newblock {\em Amer. Math. Monthly}, 73:1--23, 1966.

\bibitem{Ivrii1980:laplace_beltrami}
V.~Ja. Ivri{\u\i}.
\newblock The second term of the spectral asymptotics for a
  {L}aplace-{B}eltrami operator on manifolds with boundary.
\newblock {\em Funktsional. Anal. i Prilozhen.}, 14(2):25--34, 1980.
\newblock English translation: Functional Anal. Appl. \textbf{14} (1980), no.
  2, 98--106.

\bibitem{Berry1980:wave_motion}
M.~V. Berry.
\newblock Some geometric aspects of wave motion: wavefront dislocations,
  diffraction catastrophes, diffractals.
\newblock In {\em Geometry of the {L}aplace operator ({P}roc. {S}ympos. {P}ure
  {M}ath., {U}niv. {H}awaii, {H}onolulu, {H}awaii, 1979)}, Proc. Sympos. Pure
  Math., XXXVI, pages 13--28. Amer. Math. Soc., Providence, R.I., 1980.

\bibitem{Berry1979:distribution_modes_fractal}
M.~V. Berry.
\newblock Distribution of modes in fractal resonators.
\newblock In {\em Structural stability in physics ({P}roc. {I}nternat.
  {S}ymposia {A}ppl. {C}atastrophe {T}heory and {T}opological {C}oncepts in
  {P}hys., {I}nst. {I}nform. {S}ci., {U}niv. {T}\"ubingen, {T}\"ubingen,
  1978)}, volume~4 of {\em Springer Ser. Synergetics}, pages 51--53. Springer,
  Berlin, 1979.

\bibitem{Brossard_Carmona1986:can_one_hear}
J.~Brossard and R.~Carmona.
\newblock Can one hear the dimension of a fractal?
\newblock {\em Commun. Math. Phys.}, 104:103--122, 1986.

\bibitem{Lapidus1991:fractal_drum_inverse}
M.~L. Lapidus.
\newblock Fractal drum, inverse spectral problems for elliptic operators and a
  partial resolution of the {W}eyl-{B}erry conjecture.
\newblock {\em Trans. Amer. Math. Soc.}, 325:465--529, 1991.

\bibitem{Lapidus_Pomerance1993:riemann_weyl_berry}
M.~L. Lapidus and C.~Pomerance.
\newblock The {R}iemann zeta-function and the one-dimensional {W}eyl-{B}erry
  conjecture for fractal drums.
\newblock {\em Proc. London Math. Soc. (3)}, 66(1):41--69, 1993.

\bibitem{Lapidus_Pomerance1996:counterexamples_weyl_berry}
M.~L. Lapidus and C.~Pomerance.
\newblock Counterexamples to the modified {W}eyl-{B}erry conjecture on fractal
  drums.
\newblock {\em Math. Proc. Cambridge Philos. Soc.}, 119(1):167--178, 1996.

\bibitem{Teplyaev2007:spectral_zeta_functions}
A.~Teplyaev.
\newblock Spectral zeta functions of fractals and the complex dynamics of
  polynomials.
\newblock {\em Trans. Amer. Math. Soc.}, 359(9):4339--4358 (electronic), 2007.

\bibitem{Ransford1995:potential_theory_complex_plane}
T.~Ransford.
\newblock {\em Potential {T}heory in the {C}omplex {P}lane}, volume~28 of {\em
  London Mathematical Society Student Texts}.
\newblock Cambridge University Press, Cambridge, 1995.

\bibitem{Derfel_Grabner_Vogl2008:complex_asymptotics_julia}
G.~Derfel, P.~J. Grabner, and F.~Vogl.
\newblock Complex asymptotics of {P}oincar\'e functions and properties of
  {J}ulia sets.
\newblock {\em Math. Proc. Cambridge Philos. Soc.}, 145:699--718, 2008.

\bibitem{Hamilton1995:length_julia_curves}
D.~H. Hamilton.
\newblock Length of {J}ulia curves.
\newblock {\em Pacific J. Math.}, 169:75--93, 1995.

\bibitem{Tenenbaum1995:analytic_number_theory}
G.~Tenenbaum.
\newblock {\em Introduction to analytic and probabilistic number theory},
  volume~46 of {\em Cambridge Studies in Advanced Mathematics}.
\newblock Cambridge University Press, Cambridge, 1995.
\newblock Translated from the second French edition (1995) by C. B. Thomas.

\bibitem{He_Lapidus1997:generalized_minkowski_content}
C.~Q. He and M.~L. Lapidus.
\newblock Generalized {M}inkowski content, spectrum of fractal drums, fractal
  strings and the {R}iemann zeta-function.
\newblock {\em Mem. Amer. Math. Soc.}, 127(608):x+97, 1997.

\bibitem{Lapidus_Frankenhuysen1999:complex_dimensions_fractal}
M.~L. Lapidus and M.~van Frankenhuysen.
\newblock Complex dimensions of fractal strings and oscillatory phenomena in
  fractal geometry and arithmetic.
\newblock In {\em Spectral problems in geometry and arithmetic (Iowa City, IA,
  1997)}, volume 237 of {\em Contemp. Math.}, pages 87--105. Amer. Math. Soc.,
  Providence, RI, 1999.

\bibitem{Kumagai1993:estimates_brownian}
T.~Kumagai.
\newblock Estimates of transition densities for {B}rownian motion on nested
  fractals.
\newblock {\em Probab. Theory Related Fields}, 96(2):205--224, 1993.

\bibitem{Kajino2011:on-diagonal_oscillations}
N.~Kajino.
\newblock On-diagonal oscillation of the heat kernels on post-critically finite
  self-similar fractals.
\newblock preprint, 2011.

\bibitem{Hardy_Riesz1964:general_theory_dirichlet}
G.~H. Hardy and M.~Riesz.
\newblock {\em The general theory of {D}irichlet's series}.
\newblock Cambridge Tracts in Mathematics and Mathematical Physics, No. 18.
  Stechert-Hafner, Inc., New York, 1964.

\bibitem{Derfel_Grabner_Vogl2007:asymptotics_poincare_functions}
G.~Derfel, P.~J. Grabner, and F.~Vogl.
\newblock Asymptotics of the {P}oincar\'e functions.
\newblock In D.~Dawson, V.~Jaksic, and B.~Vainberg, editors, {\em Probability
  and Mathematical Physics: A Volume in Honor of Stanislav Molchanov},
  volume~42 of {\em CRM Proceedings and Lecture Notes}, pages 113--130,
  Montreal, 2007. Centre de Recherches Math\'ematiques.

\bibitem{Kigami_Lapidus1993:weyl's_problem_spectral}
J.~Kigami and M.~L. Lapidus.
\newblock Weyl's problem for the spectral distribution of {L}aplacians on
  p.c.f. self-similar fractals.
\newblock {\em Comm. Math. Phys.}, 158:93--125, 1993.

\bibitem{Feller1966:introduction_probability_theory_ii}
W.~Feller.
\newblock {\em An Introduction to Probability Theory and Its Applications},
  volume~II.
\newblock J. Wiley, 1966.

\bibitem{Levitin_Vassiliev1996:spectral_asymptotics}
M.~Levitin and D.~Vassiliev.
\newblock Spectral asymptotics, renewal theorem, and the {B}erry conjecture for
  a class of fractals.
\newblock {\em Proc. London Math. Soc. (3)}, 72:188--214, 1996.

\bibitem{Kajino2010:spectral_laplacians}
N.~Kajino.
\newblock Spectral asymptotics for {L}aplacians on self-similar sets.
\newblock {\em J. Funct. Anal.}, 258(4):1310--1360, 2010.

\bibitem{Poincare1886:une_classe_etendue}
H.~Poincar\'e.
\newblock Sur une classe \'etendue de transcendantes uniformes.
\newblock {\em C. R. Acad. Sci. Paris}, 103:862--864, 1886.

\bibitem{Poincare1890:une_classe_nouvelle}
H.~Poincar\'e.
\newblock Sur une classe nouvelle de transcendantes uniformes.
\newblock {\em J. Math. Pures Appl. IV. Ser.}, 6:316--365, 1890.

\bibitem{Valiron1923:lectures_on_general}
G.~Valiron.
\newblock {\em Lectures on the {G}eneral {T}heory of {I}ntegral {F}unctions}.
\newblock Private, Toulouse, 1923.

\bibitem{Valiron1954:fonctions_analytiques}
G.~Valiron.
\newblock {\em Fonctions {A}nalytiques}.
\newblock Presses Universitaires de France, Paris, 1954.

\bibitem{Eremenko_Levin1989:periodic_points_polynomials}
A.~{\`E}. Er{\"e}menko and G.~M. Levin.
\newblock Periodic points of polynomials ({R}ussian).
\newblock {\em Ukrain. Mat. Zh.}, 41:1467--1471, 1581, 1989.
\newblock translation in Ukrainian Math. J. \textbf{41} (1989), 1258--1262.

\bibitem{Eremenko_Levin1992:estimation_characteristic_exponents}
A.~{\`E}. Er{\"e}menko and G.~M. Levin.
\newblock Estimation of the characteristic exponents of a polynomial
  ({R}ussian).
\newblock {\em Teor. Funktsi\u\i{} Funktsional. Anal. i Prilozhen.}, 58:30--40
  (1993), 1992.
\newblock translation in J. Math. Sci. (New York) \textbf{85} (1997),
  2164--2171.

\bibitem{Eremenko_Sodin1990:iterations_rational_functions}
A.~{\`E}. Er{\"e}menko and M.~L. Sodin.
\newblock Iterations of rational functions and the distribution of the values
  of {P}oincar\'e functions ({R}ussian.
\newblock {\em Teor. Funktsi\u\i{} Funktsional. Anal. i Prilozhen.}, 53:18--25,
  1990.
\newblock translation in J. Soviet Math. \textbf{58} (1992), 504--509.

\bibitem{Ishizaki_Yanagihara2005:borel_and_julia}
K.~Ishizaki and N.~Yanagihara.
\newblock Borel and {J}ulia directions of meromorphic {S}chr\"oder functions.
\newblock {\em Math. Proc. Cambridge Philos. Soc.}, 139:139--147, 2005.

\bibitem{Levin1991:pommerenke's_inequality}
G.~M. Levin.
\newblock On {P}ommerenke's inequality for the eigenvalues of fixed points.
\newblock {\em Colloq. Math.}, 62:167--177, 1991.

\bibitem{Romanenko_Sharkovsky2000:long_time_properties}
E.~Romanenko and A.~Sharkovsky.
\newblock Long time properties of solutions of simplest $q$-difference
  equations ({R}ussian).
\newblock Preprint, 2000.

\bibitem{Falconer2003:fractal_geometry}
K.~J. Falconer.
\newblock {\em Fractal {G}eometry}.
\newblock John Wiley \& Sons Inc., Hoboken, NJ, 2003.
\newblock Mathematical foundations and applications.

\bibitem{Pommerenke1986:conformal_mapping_iteration}
Ch. Pommerenke.
\newblock On conformal mapping and iteration of rational functions.
\newblock {\em Complex Variables Theory Appl.}, 5:117--126, 1986.

\bibitem{Buff2003:bieberbach_conjecture_dynamics}
X.~Buff.
\newblock On the {B}ieberbach conjecture and holomorphic dynamics.
\newblock {\em Proc. Amer. Math. Soc.}, 131:755--759 (electronic), 2003.

\end{thebibliography}
